\documentclass{amsart}

\usepackage{amsthm, amssymb, amscd, a4, mathptmx}  

\DeclareSymbolFont{rsfs}{U}{rsfs}{m}{n}
\DeclareSymbolFontAlphabet{\mathscr}{rsfs}

\renewcommand{\Bbb}{\mathbb}
\renewcommand{\frak}{\mathfrak}
\newcommand{\cal}{\mathscr}
\newcommand{\catqot}{/\hskip-3pt/} \newcommand{\C}{{\Bbb C}}
\newcommand{\E}{{\cal E}}
\newcommand{\ev}{\mathop{\rm ev}}

\newcommand{\F}{{\cal F}}
\newcommand{\G}{{\cal G}}
\newcommand{\GL}{\mathop{\rm GL}}
\renewcommand{\H}{{\cal H}}
\newcommand{\gr}{\mathop{\rm gr}}
\newcommand{\Hom}{\mathop{\rm Hom}}
\newcommand{\id}{\mathop{\rm id}}
\newcommand{\K}{{\cal K}}
\renewcommand{\L}{{\cal L}}
\newcommand{\la}{\lambda}

\renewcommand{\O}{{\cal O}}
\renewcommand{\P}{{\Bbb P}}
\newcommand{\Pic}{\mathop{\rm Pic}}
\newcommand{\PGL}{\mathop{\rm PGL}}
\newcommand{\Q}{\Bbb Q}
\newcommand{\SL}{\mathop{\rm SL}}

\newcommand{\Z}{{\Bbb Z}}
\newcommand{\N}{{\Bbb N}}
\newcommand{\lra}{\longrightarrow}
\newcommand{\n}{{\cal N}}

\newcommand{\p}{\prime}
\newcommand{\q}{\quad}
\renewcommand{\phi}{\varphi}
\newcommand{\rk}{\mathop{\rm rk}}
\newcommand{\eps}{\varepsilon}

\newcommand{\Tr}{\mathop{\rm Trace}}
\theoremstyle{plain}
\newtheorem{Thm}{\q\sc Theorem}[section]
\newtheorem{Cor}[Thm]{\q\sc Corollary}
\newtheorem{Prop}[Thm]{\q\sc Proposition}
\newtheorem{Lem}[Thm]{\q\sc Lemma}
\newtheorem{Claim}{\q\sc Claim}

\newtheorem*{Thm*}{\q\sc Theorem}

\theoremstyle{remark}
\newtheorem{Rem}[Thm]{\q Remark}
\newtheorem{Ass}[Thm]{\q Assumptions}

\newtheorem{Ex}[Thm]{\q\it Example}

\title{framed hitchin pairs}

\author{Alexander H.W.~Schmitt}
\thanks{The paper is to appear in the Revue roumaine de math\'ematiques pures et appliqu\'ees}
\begin{document}
\begin{abstract}
We provide a construction of the moduli spaces of
framed Hitchin pairs and their master spaces. These objects
have come to interest as algebraic versions of solutions of
certain coupled vortex equations. Our method unifies
and generalizes constructions of several similar moduli spaces.
\\
\\
2000 \it Mathematics Subject Classification\rm: 14D20, 14D22, 14L30.
\end{abstract}
\maketitle
%\markboth{}{}
%
\section*{Introduction}
The theory of stable vector bundles on complex projective
manifolds has two completely different aspects:
The algebro-geometric part of the
construction of their moduli spaces by GIT,
on the one hand, and, on the other
hand, the theory of Hermitian-Einstein bundles which is differential
geometry. The two theories are related by the famous Kobayashi-Hitchin
correspondence. Using this, one could compute Donaldson
invariants with the help of algebraic geometry.
Now, Kobayashi-Hitchin type correspondences occur in many other
places, e.g., for Bradlow pairs, Higgs pairs, and oriented pairs.
Thus, it is desirable to have algebraic moduli spaces
for the respective objects.
For the above examples, the moduli spaces of the corresponding
stable objects in algebraic geometry have been constructed
(see, e.g., \cite{Th}, \cite{HL}, \cite{Ni}, \cite{Yo}, \cite{Sch},
\cite{OST}).
\par
In this paper,
we study
framed Hitchin pairs and oriented framed Hitchin pairs from the
algebro-geometric viewpoint. Let $X$ be a smooth projective curve and
fix line bundles $L$, $M$, and a vector bundle $H$ on $X$.
Then a framed Hitchin pair consists of a vector bundle
$E$ with $\det(E)\cong M$,
a complex number $\eps$,
a twisted endomorphism $\phi\colon E\lra E\otimes L$, and a framing
$\psi\colon E\lra H$,  and an oriented
framed Hitchin pair consists of $(E,\eps,\phi,\psi)$
as before, and an orientation $\delta\colon \det(E)\lra M$.
The motivation to study these objects comes from non-abelian
Seiberg-Witten theory as explained in the recent thesis
of M.-S.~Stupariu \cite{Stup}. He starts with the U$(2)$-
and $\mathop{\rm PU}(2)$-monopole equations on a hermitian rank two
vector bundle on a K\"ahler surface and then applies the method of
dimensional reduction to get new equations in complex dimension one.
These are certain vortex-type equations coupled with Higgs fields.
In the algebro-geometric setting, the
U$(2)$-
and $\mathop{\rm PU}(2)$-monopole equations correspond to Bradlow pairs
and oriented pairs, respectively, and --- as in Hitchin's work ---
the process of dimensional reduction has the effect of "adding" a
trace-free twisted endomorphism $\phi\colon E\lra E\otimes K_X$ (see
\cite{Stup}, p.29ff).
In the spirit of the theory of complex vector bundles,
one must now generalize the concept of Einstein-metrics to the
respective differential-geometric objects, introduce a suitable
stability
concept for (oriented) framed Hitchin pairs, and then relate both sides
by a Kobayashi-Hitchin correspondence.
For framed Hitchin pairs, this was carried out by Lin~\cite{Lin}.
The corresponding results for oriented framed Hitchin pairs are the
content of Stupariu's thesis which also contains a discussion of
Lin's results.
\par
It is therefore important to construct the algebraic models for the moduli
spaces of stable (oriented) framed Hitchin pairs together with their Gieseker
compactifications.
This will be the main concern of our note.
We will define a very general notion of framed Hitchin pairs over
arbitrary base manifolds, explain the correct notions of semistability,
and carry out a construction of the moduli spaces, using Geometric
Invariant Theory, restricting ourselves to the case
of curves in the case of oriented framed Hitchin pairs. 
In contrast to other constructions of very similar moduli spaces 
(\cite{Ni}, \cite{Yo}), for framed Hitchin pairs, 
we will follow a Simpson-type construction,
generalizing the one in \cite{Sch}.
The advantage of this approach is twofold. First, one immediately 
obtains projective moduli spaces, and second, the symmetricity
condition that the induced homomorphism $\bigwedge^2\G\otimes\E\lra\E$
be zero, usually appearing in this context, can be suppressed.
%%%
\section*{Preliminaries}
%%%
We work over $\C$. Here is a list of data which have to be fixed. We
will refer to this list without further notice.
%%%
\begin{itemize}
\item $X$, a smooth, projective scheme over the complex numbers,
\item an ample sheaf $\O_X(1)$,
\item a Hilbert polynomial $P$, $r$ and $d$ will denote the rank
      and the degree w.r.t.\ $\O_X(1)$ which this polynomial determines,
      $\mu:=d/r$,
\item a locally free sheaf $\G$ on $X$,
\item a torsion free sheaf $\H$ on $X$,
\item a Poincar\'e line bundle ${\cal N}$ on $\Pic X\times X$.
\end{itemize}
%%%%
For any coherent sheaf $\E$ on $X$, $P_\E$ and $P(\E)$ stand for its
Hilbert polynomial w.r.t.\ $\O_X(1)$.
For any scheme $S$ and any $S$-flat coherent sheaf ${\frak E}_S$ on
$S\times X$,
there exists a morphism
$d({\frak E}_S)\colon S\lra \Pic X$, associated to the
line bundle $\det({\frak E}_S)$ on $S\times X$.
We write ${\cal N}[{\frak E}_S]:=(d({\frak E}_S)\times\id_X)^*{\cal N}$.
%%%
\subsection*{\q \sc A universal construction}
%%%
Remember the following standard
construction which will be used frequently in our note:
%%%
\begin{Prop}
\label{Constr}
Let $T$ be a noetherian scheme,
${\frak E}^1_T$ and ${\frak E}_T^2$ $T$-flat coherent sheaves on
$T\times X$, and $\phi_T\colon {\frak E}^1_T\lra
{\frak E}_T^2$
a homomorphism.
Then there is a closed subscheme ${\frak V}\subset T$ whose closed points are
those $t\in T$ for which $\phi_{T|\{t\}\times X}\equiv 0$.
\end{Prop}
%%%
\subsection*{\q \sc The commutation principle}
%%%
Let us recall two results from \cite{OST} which will be used in the
second part of our paper.
%%%
\begin{Thm}
\label{Commute}
Let $G$ be a reductive group without characters,
and suppose we are given a $(\C^*\times G)$-action on the
projective scheme ${\frak X}$ which is linearized in, say, ${\frak M}$.
Then the following conditions are equivalent:
\begin{enumerate}
\item The point $x$ is $G$-semistable w.r.t.\ the given linearization.
\item There exists a linearization $l$ of the $\C^*$-action in
      some power ${\frak M}^{\otimes m}$
      such that $x$ is $l$-semistable, and the image of $x$ in
      ${\frak X}\catqot_l \C^*$ is a $G$-semistable point w.r.t.\ the
      induced
      linearization.
\end{enumerate}
\end{Thm}
%%%
\it Proof\rm. \cite{OST}, Thm.1.4.1.~\&~Rem.1.1.1.
$\square$
\par\bigskip
%%%
\begin{Rem} Note that the same statement holds for polystable points but
not for stable points.
\end{Rem}
%%%
\begin{Prop}
\label{IndPol}
Let $\C^*$ act on a vector space $W_1\oplus W_2$
with weights $e_1$ and $e_2$, $e_1<e_2$.
Then, for the linearization of the $\C^*$-action
given by $k\in\Z_{>0}$, $e\in \Z$, with
$e_1<e/k<e_2$, of the resulting $\C^*$-action on $\P(W_1^\vee\oplus
W_2^\vee)$ in $\O(k)$, the quotient is $\P(W_1^\vee)\times\P(W_2^\vee)$, and 
$\O(k(e_2-e_1))$ descends to $\O(ke_2-e,-ke_1+e)$. Futhermore, 
the linearizations with $e/k=e_i$ yield as quotients
$\P(W_i^\vee)$, $i=1,2$, $\O(k(e_2-e_1))$ descending to
$\O(k(e_2-e_1))$.
\end{Prop}
%%%
\it Proof\rm. \cite{OST}, Example~1.2.5.
$\square$\par\bigskip
%%%
%
\subsection*{\q\sc The problem of non-commuting matrices}
Denote by $M_r$ the vector space of complex $(r\times r)$-matrices,
and by $\SL_r\subset M_r$ the special linear group.
We are interested in the right action of $\SL_r$ on
$W:=M_r^{\oplus u}$ by conjugation, i.e., $(m_1,...,m_u)\cdot g=
(g^{-1}m_1 g,...,\allowbreak g^{-1} m_u g)$, for $(m_1,...,m_u)\in W$
and $g\in \SL_r$.
First of all, it is very easy to describe the (semi)stable points: 
%%%
\begin{Lem}
\label{GIT}
A point $\underline{m}:=(m_1,...,m_u)$ in $W$ is unstable
if and only if the $m_i$'s can be simultaneously triangularized,
and $\underline{m}$ fails to be semistable, if, in addition,
all the $m_i$'s are nilpotent.
\end{Lem}
%%%
\begin{Rem}
We can formulate the condition of being a nullform in another way.
For this, we think of $\underline{m}$ as a linear map from
$\C^r$ to $\C^r\otimes\C^u$. Then, $\underline{m}$ is a nullform
if and only if $(\underline{m}\otimes {\id}_{{\C^u}^{\otimes r-1}})\circ
\cdots\circ \underline{m}=0$.
When the latter happens, we call $\underline{m}$ \it nilpotent\rm.
\end{Rem}
%%%
For us, it will be important to know the ring of invariants
$\C[W]^{\SL_r}=\C[W]^{\GL_r}$.
From a technical point of view, the above problem is just a matrix problem
associated with a quiver, namely with the one consisting of one vertex
and $u$ arrows, connecting the vertex to itself.
Coordinate rings of general quiver varieties have been explicitly
determined by Le Bruyn and Procesi \cite{LP}, and Lusztig \cite{Lu}.
First, let us define some invariants. For this, let ${\Bbb F}_u$
be the free group in $u$ generators $x_1,...,x_u$.
We think of the elements of ${\Bbb F}_u$ as words in $u$ letters.
An element $\omega\in {\Bbb F}_u$ and an element
$(m_1,...,m_u)\in M_r^u$ define a matrix $m_\omega$ which is obtained
by substituting $m_i$ for the indeterminate $x_i$, $i=1,...,u$, and we 
associate to $\omega$ the invariant $T_\omega$ 
which assigns to $(m_1,...,m_u)$ the Trace of $m_\omega$.
Theorem~1 of \cite{LP} can be stated as follows in our context.
%%%
\begin{Thm}
\label{Gens}
The algebra $\C[W]^{\SL_r}$ is generated by the elements $T_\omega$
belonging to words $\omega$ of length at most $r^2$.
\end{Thm}
%%%
\begin{Rem}
As a dimension count shows, the number of invariants $s$ is much bigger
than the dimension of the quotient, whence a lot of relations must
hold. These follow all from the
Cayley-Hamilton theorem.
\end{Rem}
%%%
\begin{Ex}
[The case $r=2=u$]
%%%
For $(m_1,m_2)\in
M_2\oplus M_2=:W$, we define 
the invariants 
\begin{eqnarray*}
{\widetilde{T}}_1(m_1,m_2)&:=&\Tr(m_1),\qquad  
{\widetilde{T}}_2(m_1,m_2)\q:=\q\det(m_1),\\
{\widetilde{T}}_3(m_1,m_2)&:=&\Tr(m_2),\qquad
{\widetilde{T}}_4(m_1,m_2)\q:=\q\det(m_2),\\ \hbox{and}\q
{\widetilde{T}}_5(m_1,m_2)&:=&\Tr(m_1m_2).
\end{eqnarray*}
One can then verify that $
\C[W]^{\SL_2}= \C[{\widetilde{T}}_1,...,{\widetilde{T}}_5]$
and that this ring is isomorphic to the polynomial ring in five
variables.
\end{Ex}
%%%
Let $s$ be the number of (non-empty) words of length at most $r^2$ in
$u$ letters, $T_1,...,T_s$ be the generating invariants from Thm.~\ref{Gens}, 
and $\Xi:= {\Bbb A}^s$.
Thus, we can associate to any
$\underline{m}$ the element $\xi(\underline{m}):=
(\, T_1(\underline{m}),...,T_s(\underline{m})\,)\in \Xi$
which we call \it the characteristic vector \rm of $\underline{m}$.
This is the replacement for the
charateristic polynomials in the case of commuting matrices.
The results of this paragraph can be elegantly stated as
\par
\smallskip
{\it $\underline{m}$ is not nilpotent if and only if its
characteristic vector is non-zero.}
%%%
\section{Framed Hitchin pairs}
%%%
\subsection*{\q\sc Definitions}
We will now introduce framed Hitchin pairs.
They form a large class of objects which
comprises all the objects studied in \cite{HL},
\cite{Ni}, \cite{Yo}, and \cite{Sch}.
%%%
\subsubsection*{\q Framed Hitchin pairs}
%%%
A \it framed Hitchin pair of type
\rm
$(P,\G,\H)$ is a quadruple
$(\E,\eps,\phi,\psi)$  composed of the following ingredients
\begin{itemize}
\item a torsion free coherent sheaf $\E$ with $P_\E=P$,
\item a complex number $\eps$,
\item a twisted endomorphism $\phi\colon
\E\lra
\E\otimes
\G$,
\item a non-zero framing $\psi\colon\E\lra\H$,
\end{itemize}
such that $\phi$ is not nilpotent, i.e., $(\phi\otimes\id_{\G^{\otimes
i-1}})\circ\cdots\circ\phi\neq 0$, for all $i\in\N$, when $\eps=0$.
To make the condition of nilpotency a bit more transparent, we state
the following
%%%
\begin{Lem}
Let $\E$ be a torsion free sheaf of rank $r$ and $\phi\colon
\E\lra\E\otimes\G$ a twisted endomorphism. Then $\phi$ is nilpotent
in the above sense if and only if
$
(\phi\otimes\id_{\G^{\otimes
r-1}})\circ\cdots\circ\phi= 0$.
\end{Lem}
%%%
\it Proof\rm.
One direction is trivial. Suppose $\phi$ is nilpotent and set
$\F_i:=\ker\bigl((\phi\otimes{\id}_{\G^{\otimes
i-1}})\circ\cdots\circ\phi\bigr)$, $i\in\N$, so that
we get a filtration
$
0\subset\F_1\subset\cdots\subset\F_s=\E$.
Now, the $\F_i$'s are saturated subsheaves of $\E$. Hence, all
the inclusions above are either equalities or the rank jumps by one.
Thus, the assertion  follows from the obvious fact $\phi(\F_i)\subset
\F_{i-1}\otimes
\G$.
$\square$\par\bigskip
%%%
The equivalence relation $\sim$ on framed Hitchin pairs is
the one generated by
$$
\begin{array}{l}
(\E,\eps,\phi,\psi)\sim
(\E^\p,\eps,(\rho\otimes{\id}_\G)\circ\phi\circ\rho^{-1},
\psi\circ\rho^{-1}),\quad \rho\colon \E\lra\E^\p\ \ \hbox{iso}.,\\
(\E,\eps,\phi,\psi)\sim (\E,z\cdot\eps,z\cdot\phi,\psi),\qquad z\in\C^*.
\end{array}
$$
\begin{Rem}
Applying the above definition to the automorphism $\la\cdot \id$,
$\la\in\C^*$, one sees that $(\E,\eps,\phi,\psi)$ is always equivalent
to $(\E,\eps,\phi,\la\cdot\psi)$.
\end{Rem}
%%%
A \it family of framed Hitchin pairs of type $(P,\G,\H)$ parametrized by
the
(noetherian) scheme $S$ \rm is a quintuple $({\frak
E}_S,\eps_S,\phi_S,\psi_S,{\frak N}_S)$. Here, ${\frak E}_S$ is an
$S$-flat family of torsion free coherent sheaves on $S\times X$ with
Hilbert polynomial
$P$, ${\frak N}_S$ is a line bundle on $S$, $\eps_S\in H^0({\frak
N}_S)$, $\phi_S\colon {\frak E}_S\lra {\frak E}_S\otimes\pi_S^*{\frak
N}_S\otimes \pi_X^*\G$ is a twisted
endomorphism, and
$\psi_S\colon
{\frak E}_S\lra \pi_X^* \H$ is a framing.
We say that the family
$({\frak
E}_S,\eps_S,\phi_S,\psi_S,{\frak N}_S)$  is \it equivalent \rm to the
family
$
({\frak
E}^\p_S,\eps^\p_S,\phi^\p_S,\psi^\p_S,\allowbreak {\frak N}^\p_S)$,
if we find isomorphisms $\rho_S\colon {\frak E}_S\lra {\frak E}_S^\p$
and $\eta_S\colon {\frak N}_S\lra {\frak N}_S^\p$, such that
$$
\eps_S^\p=\eta_S\circ \eps_S,\q \phi_S^\p=(\rho_S\otimes
\pi_S^*\eta_S\otimes {\id}_{\pi_X^*\G})\circ\phi_S\circ \rho_S^{-1},\q
\psi^\p_S=\psi_S\circ\rho_S^{-1}.
$$
%%%
\subsubsection*{\q Symmetric framed Hitchin pairs and characteristic
polynomials}
%%%
Let $\E$ be a torsion free coherent sheaf and
$\phi\colon\E\lra\E\otimes\G$ be a twisted endomorphism.
We call $\phi$ \it symmetric\rm, if the induced homomorphism
$\E\otimes T^*(\G^\vee)\lra \E$ factors through
$\E\otimes S^*(\G^\vee)$.
%%%
\begin{Rem}
\label{Symm}
Since the symmetric algebra $S^*(\G^\vee)$is generated over $\O_X$ by
$\G^\vee$ it is sufficient to have that the map from
$\E\otimes\G^\vee\otimes\G^\vee\lra \E$ vanishes on $\E\otimes
(\bigwedge^2\G^\vee)$. Moreover, if $\G=\O_X(m)^{\oplus u}$
for some $m$ and $u$, we can decompose $\phi$ into its components
$(\phi_1,...,\phi_u)$. The condition of symmetricity means that the
$\phi_i$ commute, i.e., for all $i,j=1,...,u$, $(\phi_i\otimes
\id_\G)\circ \phi_j-
(\phi_j\otimes
\id_\G)\circ\phi_i=0$.
\end{Rem}
%%%
As Yokogawa explains in \cite{Yo2}, p.495, we can associate
to a symmetric twisted endomorphism $\phi$ its characteristic
polynomial in $H^0(S^*(\G)[t])$. This provides us, in particular, with
an element in the vector space ${\Bbb H}_\G:=\bigoplus_{i=1}^r H^0(S^i(\G))$ which we
call abusively the \it characteristic polynomial of $\phi$\rm , too.
Note that $\phi$ is not nilpotent if and only if its characteristic
polynomial in ${\Bbb H}_\G$ is not zero.
\par
A framed Hitchin pair $(\E,\eps,\phi,\psi)$ of type $(P,\G,\H)$
will be called \it symmetric\rm , if $\phi$ is symmetric.
A symmetric framed Hitchin pair $(\E,\eps,\phi,\psi)$ of type
$(P,\G,\H)$ defines an element in $\C\oplus {\Bbb H}_\G$.
Now, let $\C^*$ act on $H^0(S^i(\G))$ through multiplication by $z^i$.
This yields a $\C^*$-action on $\C\oplus {\Bbb H}_\G$, and the quotient
$\widehat{\P}_\G$ is a weighted projective space.
Thus, $(\E,\eps,\phi,\psi)$ defines an element
$\widehat{\chi}(\E,\eps,\phi,\psi)\in\widehat{\P}_\G$ which depends only
on
its equivalence class. It will be referred to as \it the characteristic polynomial
of
$(\E,\eps,\phi,\psi)$\rm.
%%%
%%%
\subsubsection*{\q Characteristic vectors for framed Hitchin pairs}
%%%
Let $(\E,\eps,\phi,\psi)$ be a framed Hit\-chin pair of type $(P,\G,\H)$
with $\G=\O_X(m)^{\oplus u}$. Let $T_1,...,T_s$ be the generators of the 
ring $\C[M_r^{\oplus u}]^{\SL_r}$, belonging to the words
$\omega_1,....,\omega_s$.
For a given $l$, consider the homomorphism $\phi_l:=(\phi\otimes 
{\id}_{\G^{\otimes l-1}})\circ \cdots \circ \phi\colon \E\lra \E\otimes 
\G^{\otimes l}$. Now, a word $\omega$ of length $l$ singles out a 
component of $\phi_l$, i.e., a homomorphism $\phi_\omega\colon \E\lra 
\E\otimes \O_X(lm)$, and the trace of $\phi_\omega$ gives a section
in $H^0(\O_X(lm))$. For $i=1,...,s$, let $l_s$ be the length of the word
$\omega_s$, and set ${\Bbb K}_\G:=\bigoplus_{i=1}^s H^0(\O_X(l_im))$.
So, we can associate to $(\E,\eps,\phi,\psi)$ in a natural way
an element in ${\Bbb K}_\G$.
As before, we let $\C^*$ act on $H^0(\O_X(l_im))$ through multiplication
by $z^{l_i}$, $i=1,...,s$, in order to obtain a $\C^*$-action
on $\C\oplus {\Bbb K}_\G$. The quotient $(\C\oplus {\Bbb K}_G)\catqot\C^*$
will be denoted by $\widehat{\Xi}_\G$.
Thus, we can assign to any framed Hitchin pair $(\E,\eps,\phi,\psi)$
its \it characteristic vector \rm $\widehat{\xi}(\E,\eps,\phi,\psi)\in
\widehat{\Xi}_\G$.
The characteristic vector clearly depends only on the equivalence class
of $(\E,\eps,\phi,\psi)$.
%%%
\subsection*{\q\sc Semistability and sectional semistability}
%%%
Fix a polynomial $\sigma\in{\Q}[t]$ with positive leading coefficient
of degree at most
$\dim X-1$ and call a framed Hitchin pair $(\E,\eps,\phi,\psi)$
\it $\sigma$-(semi)stable\rm, if
for any proper, non-zero subsheaf $\F$ of $\E$ which is invariant
under $\phi$, i.e., $\phi(\F)\subset \F\otimes \G$,
\begin{eqnarray*}
{P_\F\over\rk\F}-{\sigma\over\rk\F}& (\le)&
{P_\E\over\rk\E}-{\sigma\over\rk\E},\q\hbox{and} \\
{P_\F\over\rk\F}& (\le)&
{P_\E\over\rk\E}-{\sigma\over\rk\E},  \q \hbox{if, furthermore, }
\F\subset\ker\psi.
\end{eqnarray*}
%%%
\begin{Rem}
\label{Replace}
i) Given
two locally free sheaves $\G\subset\G^\p$  and a framed Hit\-chin pair
$(\E,\eps,\phi,\allowbreak \psi)$ of type $(P,\G,\H)$, we can form the framed
Hitchin pair $(\E,\eps,\phi^\p,\psi)$ of type $(P,\G^\p,\H)$ by
defining $\phi^\p$ as the composition of $\phi$ with the inclusion
$\E\otimes\G\subset\E\otimes\G^\p$.
Since for any $\phi^\p$-invariant subsheaf $\F$ of $\E$, the map
$\F\lra\F\otimes\G^\p\lra \F\otimes (\G^\p/\G)$ is zero, it will also
be $\phi$-invariant.
Therefore, $(\E,\eps,\phi^\p,\psi)$ will be $\sigma$-(semi)stable if and only
if $(\E,\eps,\phi,\psi)$ is $\sigma$-(semi)stable.
In particular, since we can embed $\G$ in ${\O_X(m)}^{\oplus u}$
for some large $m$ and $u$, we can and will always assume that $\G$ is
of that simple form.
\par
ii)
A $\sigma$-stable framed Hitchin pair $(\E,\eps,\phi,\psi)$ 
has no automorphisms $\rho$ besides the identity
which satisfy
$$
(\E,\eps,\phi,\psi)\q =\q \bigl(\E,\eps,
(\rho\otimes{\id}_{\G})\circ\phi\circ\rho^{-1},\psi\circ\rho^{-1}\bigr).
$$
\end{Rem}
%%%
Next, fix a positive rational number $\overline{\sigma}$. A framed
Hitchin pair $(\E,\eps,\phi,\psi)$ of type $(P,\G,\H)$ will be called
\it $\overline{\sigma}$-sectional (semi)stable\rm, if there is a
subspace $V\subset H^0(\E)$ of dimension $\chi(\E)$, s.~th.\ for any
non-trivial $\phi$-invariant proper subsheaf $\F$ of $\E$
%%%
\begin{eqnarray*}
{\dim(V\cap H^0(\F))\over\rk\F} -{\overline{\sigma}\over\rk\F} &(\le)&
{\chi(\E)\over\rk\E} -{\overline{\sigma}\over\rk\E}, \\
{\dim(V\cap H^0(\F))\over\rk\F}
&(\le)&
{\chi(\E)\over\rk\E} -{\overline{\sigma}\over\rk\E},\qquad \hbox{if}
\q\F\subset\ker\psi.
\end{eqnarray*}
%%%
The usual
arguments --- assuming boundedness  --- then
show:
%%%
\begin{Prop}
\label{Sect}
There is a natural number $n_0$ such that for all $n\ge n_0$ and all
framed Hitchin pairs $(\E,\eps,\phi,\psi)$ of type $(P,\G,\H)$ the
following conditions are equivalent:
\begin{enumerate}
\item $(\E,\eps,\phi,\psi)$ is $\sigma$-(semi)stable.
\item $(\E(n),\eps,\phi\otimes\id_{\O_X(n)},\psi
\otimes\id_{\O_X(n)})$ is $\sigma(n)$-sectional (semi)stable.
\item $(\E(n),\eps,\phi\otimes\id_{\O_X(n)},\psi
\otimes\id_{\O_X(n)})$ 
satisfies the condition of  $\sigma(n)$-sec\-tional (semi)\allowbreak stability
for globally generated subsheaves.
\end{enumerate}
\end{Prop}
%%%
\begin{Rem}
If $X$ is curve and some positive rational number
$\sigma_\infty$, then --- as will follow from the results of
the section about boundedness in the second chapter
 --- one can choose $n_0$ such that the
above proposition holds true for all positive 
rational numbers $\sigma<\sigma_\infty$.
\end{Rem} 
%%%
\subsection*{\q\sc Jordan-H\"older Filtrations and S-equivalence}
%%%
Let $(\E,\allowbreak\eps,\phi,\psi)$ be a $\sigma$-semi\-stable framed Hitchin pair
of type $(P,\G,\H)$ which is not $\sigma$-stable.
Let $\E_1$ be a proper, $\phi$-invariant, destabilizing subsheaf which
is maximal
w.r.t.\ inclusion.  The homomorphisms $\phi$ and $\psi$ can be
restricted to ${\E}_1$, and
$(\E_1,\eps,\allowbreak \phi_{|{\E}_1},
\psi_{|\E_1})$ will be again $\sigma$-semistable.
If it is not $\sigma$-stable, pick a maximal destabilizing
$\phi$-invariant subsheaf ${\E_2}$ and so on in
order to get a so called \it Jordan-H\"older filtration\rm
$$
0=:\E_{m+1}\subset \E_m   \subset\cdots\subset\E_1\subset\E_0:=\E.
$$
Then, we can define \it the associated graded object \rm
$$
\gr(\E,\eps,\phi,\psi)\q:=\q\bigoplus_{i=1}^{m+1}
\bigl(\E_{i-1}/\E_{i},\eps,
\overline{\phi}_{i},\overline{\psi}_i\bigr).
$$
As usual, this is well-defined up to equivalence.
Two $\sigma$-semistable framed Hitchin pairs are called
\it S-equivalent \rm if their associated graded objects are equivalent,
and a framed Hitchin pair is said to be \it $\sigma$-polystable\rm,
if it is equivalent to its associated graded object.
%%%
\subsection*{\q\sc The main result}
Let $\mathop{\rm FH}_{P/\G/\H}^{\sigma-(s)s}$
($\mathop{\rm FH}_{P/\G/\H/\hbox{\rm symm}}^{\sigma-(s)s}$)
be the functors 
assigning to each noetherian scheme the set of equivalence classes of
families of $\sigma$-(semi)stable (symmetric) framed Hitchin pairs of
type
$(P,\G,\H)$.
\begin{Thm}
\label{Main}
{\rm i)} There are a quasi-projective scheme ${\cal
{\cal F}{\cal H}}_{P/\G/\H}^{\sigma-ss}$
and a natural transformation
$\vartheta$ of $\mathop{\rm FH}_{P/\G/\H}^{\sigma-ss}$
into the functor of points of $
{\cal
{\cal F}{\cal H}}_{P/\G/\H}^{\sigma-ss}$, s.~th.\ for any scheme
${\cal M}$ and any natural transformation $\vartheta^\p\colon
\mathop{\rm FH}_{P/\G/\H}^{\sigma-ss}\lra h_{\cal M}$, there
is a unique morphism $\tau\colon
{\cal
{\cal F}{\cal H}}_{P/\G/\H}^{\sigma-ss} \lra {\cal M}$ with
$\vartheta^\p=h(\tau)\circ\vartheta$.
\par
{\rm ii)} The space
${\cal {\cal F}{\cal H}}_{P/\G/\H}^{\sigma-ss}$  contains
an open subscheme
${\cal
{\cal F}{\cal H}}_{P/\G/\H}^{\sigma-s}$
which is a fine moduli space
for the functor
$\mathop{\rm FH}_{P/\G/\H}^{\sigma-s}$.
\par
{\rm iii)} The closed points of
${\cal
{\cal F}{\cal H}}_{P/\G/\H}^{\sigma-ss}$
 naturally correspond to the set of
S-equivalence
classes of $\sigma$-semistable framed Hitchin pairs of type $(P,\G,\H)$.
\par
{\rm iv)}
There is a proper morphism $\widehat{\xi}\colon
{\cal
{\cal F}{\cal H}}_{P/\G/\H}^{\sigma-ss}
\lra \widehat{\Xi}_\G $ --- called the {\rm generalized Hitchin map}
--- which maps a $\sigma$-polystable framed Hitchin
pair of type $(P,\G,\H)$ to its characteristic vector.
In particular, ${\cal
{\cal F}{\cal H}}_{P/\G/\H}^{\sigma-ss}$ is a projective scheme.
\par
{\rm v)} There is a closed subscheme $
{\cal {\cal F}{\cal H}}_{P/\G/\H/\hbox{\rm symm}}^{\sigma-ss}$
of  ${\cal {\cal F}{\cal H}}_{P/\G/\H}^{\sigma-ss}$, such that
the schemes 
$
{\cal {\cal F}{\cal H}}_{P/\G/\H/\hbox{\rm symm}}^{\sigma-ss}$
and
${\cal {\cal F}{\cal H}}_{P/\G/\H/\hbox{\rm symm}}^{\sigma-s}:=
{\cal {\cal F}{\cal H}}_{P/\G/\H/\hbox{\rm symm}}^{\sigma-ss}\cap
{\cal {\cal F}{\cal H}}_{P/\G/\H}^{\sigma-s}$ enjoy the analogous properties
to {\rm i)} - {\rm iii)} w.r.t.~functors
$\mathop{\rm FH}_{P/\G/\H/\hbox{\rm symm}}^{\sigma-(s)s}$.
Moreover, there is a proper morphism $\widehat{\chi}\colon
{\cal
{\cal F}{\cal H}}_{P/\G/\H/\hbox{\rm symm}}^{\sigma-ss}
\lra \widehat{\P}_\G $, the {\rm Hitchin map},
which maps a $\sigma$-polystable symmetric framed Hitchin
pair of type $(P,\G,\H)$ to its characteristic polynomial.
\end{Thm}
%%%
\subsection*{\q\sc Proof of Theorem~\ref{Main} with GIT}
\label{PF}
We will follow the usual pattern of a GIT construction.
%%%
\subsubsection*{\q Boundedness}
%%%
For a torsion free coherent sheaf $\E$, let $0=\E_0\subset\E_1
\subset\cdots\subset\E_l=\E$ be its (slope) Harder-Narasimhan
filtration. Set $\mu_{\max}(\E):=\mu(\E_1)$ and $\mu_{\min}(\E):=
\mu(\E/\E_{l-1})$.
The proof of Nitsure \cite{Ni}, Proposition~3.2, can be easily extended
to give the following
%%%
\begin{Prop}
\label{bound2}
Let $(\E,\eps,\phi,\psi)$ be a framed Hitchin pair of type
$(P,\allowbreak\G,\H)$ with $\G=\O_X(m)^{\oplus u}$, such that there is a constant
$C\ge
0$ such that
$\mu(\F)\le \mu(\E)+[C(r-\rk\F)/(r\rk\F)]$ for any
non-trivial $\phi$-invariant subsheaf $\F$ of $\E$.
Then
$$
\mu_{\max}(\E)\q\le\q \max\Bigl\{\,\mu+C, \mu+ {(r-1)^2\over
r}\deg\O_X(m)\,\Bigr\}.
$$
\end{Prop}
%%%
\begin{Rem}
\label{PolVersion}
As explained in \cite{Sch}, p.111, one can formulate a
``polynomial'' analogon to Proposition~\ref{bound2}, namely, there
is  a constant $C^\p$ such that
for any subsheaf $\F$ of $\E$ one has
$P(\F)/\rk\F\ <P(\E)/\rk\E+C^\p x^{\dim X-1}$.
\end{Rem}
%%%
Since every $\sigma$-semistable framed Hitchin pair of type $(P,\G,\H)$
satisfies the assumption of the above proposition with $C=$ leading
coefficient of $\sigma$, Maruyama's boundedness result \cite{Ma}
yields
%%%
\begin{Cor}
\label{bound3}
The isomorphy classes of torsion free coherent shea\-ves occuring
in $\sigma$-semistable framed Hitchin pairs of type $(P,\G,\H)$
form a bounded family.
\end{Cor}
%%%
Let $(\E,\eps,\phi,\psi)$ be a framed Hitchin pair of type $(P,\G,\H)$
and set $\F_i:=\ker\bigl((\phi\otimes{\id}_{\G^{\otimes
i-1}})\circ\cdots\circ\phi\bigr)$, $i=1,...,r$.
%%%
\begin{Lem}
\label{bound4}
The set of isomorphy classes of $\F_i$'s coming from
$\sigma$-semi\-stable
framed Hitchin pairs of type $(P,\G,\H)$ is also bounded.
\end{Lem}
%%%
\it Proof\rm.
The proof of this lemma will be given below.
$\square$\par\bigskip
%%%
\subsubsection*{\q Some assumptions}
As observed in Remark~\ref{Replace}, we can assume that $\G$ is of the
form $\O_X(m)^{\oplus u}$ where $\O_X(m)$ is globally generated, and,
by Corollary~\ref{bound3} and Lemma~\ref{bound4}, the following 
can be required.
%%%
\begin{Ass}
\label{Ass1}
Let $n_1$ be a natural number
such that for all $n\ge n_1$ and every $\sigma$-semistable framed
Hitchin pair $(\E,\eps,\phi,\psi)$ of type $(P,\G,\H)$
the following holds
\begin{itemize}
\item
$\H(n)$ is globally generated and without higher cohomology.
\item
$\E(n)$ is globally generated and without higher cohomology.
\item 
For $i=1,...,r$, the sheaf $\F_i(n)$ is globally generated and without
higher cohomology.
\end{itemize}
\end{Ass}
%%%
An  additional hypothesis on $n_1$ will be explained later.
%%%
Now, we have to make a "logical loop". Indeed, the last assumption
makes use of~\ref{bound4} which we have not yet proved.
Therefore, we won't use this assumption in the following construction,
prove Lemma~\ref{bound4}, and then re-enter at the beginning of this
section.
%%%
Suppose also that $n_1$ is greater than the constant $n_0$ in
Proposition~\ref{Sect}.
We may clearly assume that $n_1=0$.
%%%
\subsubsection*{\q The parameter space}
\label{ParSpace}
%%%
Let $V$ be a complex vector space
of dimension $p:=P(0)$. By Assumption~\ref{Ass1}, every torsion free
sheaf $\E$ occuring in a $\sigma$-semistable framed Hitchin
pair
of type $(P,\G,\H)$ can be written as a quotient
$q\colon V\otimes \O_X\lra\E$ where $H^0(q)$ is an isomorphism.
These quotients are parametrized by a quasi-projective scheme ${\frak
Q}_0$ which is an open subscheme of
${\frak Q}$, the projective quot scheme                    of
quotients of $V\otimes\O_X$ which have Hilbert polynomial $P$.
Let
$
{\frak q}_{\frak Q}\colon V\otimes \O_{{\frak Q}\times X}\lra
{\frak E}_{\frak Q}
$
be the universal quotient on ${\frak Q}\times X$.
We can choose a $\nu_0$ meeting the following
requirements:
%%%
\begin{Ass}
\label{Ass2} For any $\nu\ge\nu_0$, any subspace $U$ of $V$,
any $[q\colon V\otimes\O_X\lra\E]\in {\frak Q}$,
and $\E_U:=q(U\otimes\O_X)$:
\begin{itemize}
\item $\O_X(\nu)$ is globally generated and without higher cohomology.
\item The map $H^0(\O_X(\nu))\otimes H^0(\O_X(m))\lra H^0(\O_X(\nu+m))$
      is surjective.
\item $H^i(\E_U(\nu))=0$, $i>0$, and $U\otimes\O_X(\nu)\lra
H^0(\E_U(\nu))$ is
      surjective.
\end{itemize}
\end{Ass}
%%%
Hence,
$\pi_{{\frak Q}*}({\frak E}\otimes\pi_X^*\O_X(\nu_0))$ and
$\pi_{{\frak Q}*}({\frak E}\otimes\pi_X^*\O_X(\nu_0+m))$ are locally
free.
Define
$$
\widehat{\frak P}\q:=\q\P\bigl( \O_{\frak Q}\oplus
{\pi_{{\frak Q}*}({\frak E}\otimes\pi_X^*\O_X(\nu_0+m))^{\oplus
u}}^\vee\otimes \pi_{{\frak Q}*}({\frak
E}\otimes\pi_X^*\O_X(\nu_0))
\bigr).
$$
There is a tautological line bundle ${\frak N}_{\widehat{\frak P}}$
on ${\widehat{\frak P}}$, and
the tautological surjection
provides us, on $\widehat{\frak P}\times X$, with a homomorphism
$$
V\otimes \pi_X^*\O_X(\nu_0)\lra ({\frak E}_{\widehat{\frak P}}\otimes
\pi_{\widehat{\frak P}}^*{\frak
N}_{\widehat{\frak P}}\otimes \pi_X^*\O_X(\nu_0+m))^{\oplus u}
\otimes\pi_{\widehat{\frak P}}^*{\frak N}_{\widehat{\frak P}}.
$$
Here, ${\frak q}_{\widehat{\frak P}}\colon V\otimes\O_{\widehat{\frak
P}\times X}
\lra{\frak E}_{\widehat{\frak P}}$ is the pullback of the universal
quotient.
Define ${\frak P}$ as
the closed subscheme of $\widehat{\frak P}$ where this homomorphism
factorizes through ${\frak E}_{\widehat{\frak P}}\otimes \pi_X^*\O_X(\nu_0)$.
Let
$$\phi_{{\frak P}}\colon {\frak E}_{{\frak P}}\lra
({\frak E}_{{\frak P}}
\otimes\pi_X^*\O_X(m))^{\oplus u}
\otimes\pi_{{\frak P}}^*{\frak N}_{{\frak P}}
$$
be the induced homomorphism
and
${\frak N}_{\frak P}$ the restriction of ${\frak
N}_{\widehat{\frak P}}$ to ${\frak P}$.
By Assumption~\ref{Ass1}, for any $\sigma$-(semi)stable
framed Hitchin pair $(\E,\eps,\phi,\psi)$, the homomorphism
$\psi\colon \E\lra \H$ is determined by the homomorphism
$H^0(\phi)\colon \allowbreak H^0(\E)\lra H^0(\H)$.
Define the space ${\bf R}:=\P(\Hom(V,H^0(\H))^\vee)$. It is now clear how to
construct
the parameter space ${\frak R}$ as a closed subscheme of ${\frak
P}\times
{\bf R}$, and that, on ${\frak R}\times X$, there is a universal family
$({\frak E}_{\frak R},\eps_{\frak R},\phi_{\frak R},\psi_{\frak
R},{\frak N}_{\frak R})$, such that any family of $\sigma$-semistable
framed Hitchin pairs can locally be obtained as the pullback of this
universal family.  ${\frak R}_0$ will be the open subscheme sitting
over ${\frak Q}_0$, and ${\frak R}_0^{\sigma-(s)s}$ will be the  ---
a posteriori --- open subscheme parametrizing $\sigma$-(semi)stable
framed Hitchin pairs.
%%%
\subsubsection*{\q Ample line bundles on $\widehat{\frak P}$}
%%%
For the moment, write ${\frak W}_{\frak Q}$ for the vector bundle
$\pi_{{\frak Q}*}({\frak E}_{{\frak Q}}\otimes 
\pi_X^*\O_X(\nu_0+m))$.
Note ${\frak W}_{\frak Q}^\vee\cong \bigwedge^{P(\nu_0+m)-1}
{\frak W}_{\frak Q}\otimes \det{\frak W}_{\frak Q}^\vee$.
In particular, there is a surjection
$$
\bigwedge^{P(\nu_0+m)-1}(V\otimes N\otimes M)\otimes 
\det{\frak W}_{\frak Q}^\vee
\lra 
{\frak W}_{\frak Q}^\vee.
$$
The line bundle ${\frak L}_{\frak Q}:=\det{\frak W}_{\frak Q}$
on ${\frak Q}$ is very ample.
Therefore, we see that $\O_{\widehat{\frak P}}(a_1,a_2)
:=\pi^*{\frak L}_{\frak Q}^{\otimes a_1}
\otimes {\frak N}_{\widehat{\frak P}}^{\otimes a_2}$
is globally generated for $a_1\ge a_2>0$ and very ample
for $a_1>a_2>0$.
%%%
\subsubsection*{\q The group actions}
%%%
There are natural right actions of $\SL(V)$ on the schemes
${\frak Q}$ and ${\bf R}$, and the locally free sheaves
$\pi_{{\frak Q}*}({\frak E}_{\frak Q}\otimes\pi_X^*\O_X(\nu_0))$
and
$\pi_{{\frak Q}*}({\frak E}_{\frak Q}\otimes\pi_X^*\O_X(\nu_0+m))$
are naturally linearized w.r.t.\ the group action on ${\frak Q}$.
Thus, there is a natural $\SL(V)$-action from the right on
$\widehat{\frak P}\times{\bf R}$. This leaves the schemes
${\frak R}$ and ${\frak R}_0$ invariant. Moreover, the equivalence
relation on the closed points of ${\frak R}_0$ induced by this group
action is just the relation $\sim$ on framed Hitchin pairs. More
precisely,
%%%
\begin{Prop}
\label{GlueTog}
For any noetherian scheme $S$ and any two morphisms $\beta_i\colon S\lra
{\frak R}_0$, $i=1,2$, such that the pullbacks of the universal family
via $(\beta_1\times\id_X)$ and $(\beta_2\times\id_X)$ are equivalent,
there exist an \'etale covering $\tau\colon T\lra S$ and a morphism
$\Gamma\colon T\lra \SL(V)$ with $\beta_1\circ
\tau=(\beta_2\circ\tau)\cdot\Gamma$.
\end{Prop}
%%%
\begin{Rem}
\label{Lin}
Using the universal properties, it is easy to see that the
universal family on ${\frak R}\times X$ can be equipped with an
$\SL(V)$-linearization.
Now, restrict everything to ${\frak R}_0^{\sigma-s}$.
Then, by Remark~\ref{Replace}~ii), the $\SL(V)$-lineariza\-tion
induces a $\PGL(V)$-linearization of the universal
family. Since all the $\PGL(V)$-stabilizers are trivial,
the universal family descends to the quotient ${\frak R}_0^{\sigma-s}\catqot\SL(V)
= {\frak R}_0^{\sigma-s}\catqot\PGL(V)$ provided the latter
space exists. For the necessary descent theory, the reader
is referred to \cite{HL2}, p.87.
\end{Rem}
%%%
\subsubsection*{\q The semistable points
in the closure of 
${\frak R}_0$ and the proof of the main theorem}
%%%
Fix the polarization $\O(2,1,a)$ (compare \cite{HL}, p.309) on
${\frak R}$ with 
$$
{a\over 2}:= \bigl(P(\nu_0+m)-\sigma(\nu_0+m)\bigr)
{\sigma(0)\over p-\sigma(0)}-\sigma(\nu_0+m).
$$
The group action will be linearized in that line bundle.
%%%
\begin{Thm}
\label{SemStab}
Let $r:=([q\colon V\otimes\O_X\lra\E],[\eps,\phi],[\psi])$
be a point in the parameter space ${\frak R}_0$. 
Then, $r$ is (semi/poly)stable in ${\frak R}$ (w.r.t.\ the chosen
linearization) if and only if $(\E,\eps,\phi,\psi)$ is
a $\sigma$-(semi/poly)stable framed Hitchin pair of type
$(P,\G,\H)$.
Moreover, if $X$ is a curve, and $r$ is a point in the closure
of the parameter space ${\frak R}_0$, then $r$ is (semi/poly)stable in 
${\frak R}$ if and only if $r$ lies in ${\frak R}_0$ and
$(\E,\eps,\phi,\psi)$ is
a $\sigma$-(semi/poly)stable framed Hitchin pair of type
$(P,\G,\H)$.
\end{Thm}
%%%
\it Proof\rm. This theorem will be proved below.
$\square$\par\bigskip
%%%
If $X$ is a curve, then Theorem~\ref{SemStab} shows that
${\frak R}_0^{\sigma-(s)s}$ is exactly the set of
(semi)stable points in the closure of ${\frak R}_0$.
Thus, the good quotients ${\frak R}_0^{\sigma-(s)s}\catqot\SL(V)$
do exist, the space ${\frak R}_0^{\sigma-ss}\catqot\SL(V)$ is a projective scheme,
and ${\frak R}_0^{\sigma-s}\catqot\SL(V)$ is a geometric quotient.
In case $X$ is higher dimensional, Theorem~\ref{SemStab} shows
that ${\frak R}_0^{\sigma-(s)s}$ is a saturated open subet
of the semistable points in $\overline{{\frak R}_0}$, i.e.,
the closure in $\overline{{\frak R}_0}^{ss}$ of the orbit of a point 
in ${\frak R}_0^{ss}$ still lies in ${\frak R}_0^{ss}$.
Therefore, the good (geometric) quotient $
{\frak R}_0^{\sigma-(s)s}\catqot\SL(V)$ exists as a quasi-projective scheme.
Defining ${\cal
{\cal F}{\cal H}}^{\sigma-(s)s}_{P/\G/\H}:=
{\frak R}_0^{\sigma-(s)s}\catqot\SL(V)$ gives our moduli spaces.
The first assertion of Theorem~\ref{Main} is then --- as usual ---
a direct consequence of the local universal property of ${\frak R}$,
Proposition~\ref{GlueTog}, and the universal property
of the
categorical quotient. As we have already remarked in~\ref{Lin},
the universal family on ${\frak R}_0^{\sigma-s}\times X$ descends
to ${\cal {\cal F}{\cal H}}_{P/\G/\H}^{\sigma-s}\times X$, whence
${\cal {\cal F}{\cal H}}_{P/\G/\H}^{\sigma-s}$
is a fine moduli
space.
The identification of the closed points follows from the assertion
about the polystable points in Theorem~\ref{SemStab}.
Thus, i) - iii) in~\ref{Main} are settled.
For point v), we define ${\cal {\cal F}{\cal H}}^{\sigma-ss}_{P/\G/\H/\hbox{\rm
symm}}$ as the image of ${\frak R}_{0,\hbox{\rm symm}}^{\sigma-ss}$
in ${\cal {\cal F}{\cal H}}^{\sigma-ss}_{P/\G/\H}$. This space clearly has the
desired properties and coincides with
${\frak R}_{0,\hbox{\rm symm}}^{\sigma-ss}\catqot \SL(V)$.
To define the generalized Hitchin map, let $\widehat{\frak H}$ be the
geometric vector bundle associated to the
locally free sheaf
$$
\pi_{{\frak Q}*}({\frak
E}\otimes\pi_X^*\O_X(\nu_0))^\vee \otimes
{\pi_{{\frak Q}*}({\frak E}\otimes\pi_X^*\O_X(\nu_0+m))^{\oplus
u}}.
$$
The induced map
$
\bigl([\C\times\widehat{\frak H}]\setminus\{\,\hbox{zero section}\,\}\bigr)
\times {\bf R}\lra \widehat{\frak P}\times{\bf R}
$
is a good $\C^*$-quotient.
Let ${\frak K}^{\sigma-ss}_{0}$ be the preimage
of ${\frak R}^{\sigma-ss}_{0}$ under the above map.
Then, ${\cal {\cal F}{\cal H}}_{P/\G/\H}^{\sigma-ss}=
{\frak K}^{\sigma-ss}_{0}\catqot (\C^*\times\SL(V))$
is also a good quotient.
Copying the construction of \cite{Yo2}, p.496, we obtain
a morphism
$
{\frak K}^{\sigma-ss}_{0}\lra \C\times {\Bbb K}_\G
\lra \widehat{\Xi}_\G
$
which is $(\C^*\times \SL(V))$-invariant and, thus, descends
to a morphism
$$
\widehat{\xi}\colon
{\cal {\cal F}{\cal H}}_{P/\G/\H}^{\sigma-ss}
\lra \widehat{\Xi}_\G.
$$
This is the generalized Hitchin map.
To see that it is proper, let $(C,0)$ be the spectrum of a discrete
valuation ring $R$ with field of fractions $K$.
Suppose there is a map $C\lra \widehat{\Xi}_\G$ which lifts
via $\widehat{\xi}$ over $C\setminus 0$ to $
{\cal {\cal F}{\cal H}}_{P/\G/\H}^{\sigma-ss}$.
Since there are no non-trivial line bundles on $C$,
the morphism from $C$ to $\widehat{\Xi}_\G$ lifts to $\C\oplus{\Bbb
K}_\G$. After possibly passing to a finite extension of $K$,
we may assume that the map from $C\setminus 0\lra
{\cal {\cal F}{\cal H}}_{P/\G/\H}^{\sigma-ss}$ comes from
a family of $\sigma$-semistable framed Hitchin pairs of type $(P,\G,\H)$
and that the induced map to ${\Bbb K}_\G$ is just the characteristic
vector of that family. This follows from Luna's \'etale slice
theorem and our definition of equivalence of families.
It is not hard to see that the arguments used by Yokogawa 
(\cite{Yo2}, p.487ff) or Nitsure \cite{Ni} can be adapted to our
situation. We omit this here, because it does not involve any new idea.
$\square$
%%%
\subsubsection*{\q Proof of Theorem~\ref{SemStab}}
%%%
First, we prove the statement about the (semi)stable points,
using the Hilbert-Mumford criterion. For this, we recall that
a one parameter subgroup of $\SL(V)$ is determined by giving
a basis $v_1,...,v_p$ of $V$ and weights $\gamma_1\le \cdots\le
\gamma_p$, satisfying $\sum_i\gamma_i=0$.
Moreover, a weight vector $(\gamma_1,...,\gamma_p)$ with
$\gamma_1\le \cdots\le
\gamma_p$ and $\sum_i\gamma_i=0$ can be written as
$$
(\gamma_1,...,\gamma_p)\q =\q
\sum_{i=1}^{p-1}{\gamma_{i+1}-\gamma_i\over p}\gamma^{(i)}
$$
with $\gamma^{(i)}=(i-p,...,i-p,i,...,i)$ where $i-p$ occurs $i$ times.
Let us first make some preliminary remarks.
%%%
\subsubsection*{\q The weights of the ${\bf R}$-component}
%%%
Let $\la$ be the one parameter subgroup which is
determined by the basis $v_1,...,v_p$ and the weight vector
$(\gamma_1,...,\gamma_p)$. Let $[h]\in {\bf R}$ be a point.
It follows that $\mu_{\bf R}([h],\la)=-\min\{\, \gamma_i\,|\,
h(v_i)\neq 0\,\}$. In particular, if $\la^{(i)}$ is the one parameter
subgroup defined by $\gamma^{(i)}$, then $\mu_{\bf R}([h],\la^{(i)})=
-i$ or $p-i$, depending on whether $\langle\,
v_1,...,v_i\,\rangle$ is contained in $\ker(h)$ or not.
Observe that $\mu_{\bf R}([h],\la_1\cdot \la_2)=
\mu_{\bf R}([h],\la_1)+\mu_{\bf R}([h],\la_2)$ for any two
one parameter subgroups given w.r.t.\ above basis by the weights
$\gamma_1^i\le\cdots\le\gamma_p^i$,
$i=1,2$.
%%%
\subsubsection*{\q The weights of ${\frak L}_{\frak Q}$}
%%%
Let $v_1,...,v_p$ be a basis for $V$ and $[q]\in {\frak Q}$.
Define $Q_i:=H^0(q\otimes\id_{\O_X(\nu_0+m)})(\langle\,
v_1,...,v_i\,\rangle\otimes N\otimes H^0(\O_X(m)))$, 
and $\delta(i):=\dim Q_i$, $i=1,...,p$.
For a one parameter subgroup given w.r.t.\ the above basis by
the weights $\gamma_1\le\cdots\le\gamma_p$, we obtain
$$
\mu_{\frak Q}([q],\la)\q =\q
-\sum_{i=1}^p\bigl(\delta(i)-\delta(i-1)\bigr)\gamma_i,
$$
in particular, $\mu_{\frak
Q}([q],\la^{(i)})=\bigl(p\delta(i)-iP(\nu_0+m)\bigr)$. Again,
$\mu_{\frak Q}([q],\la_1\cdot \la_2)=
\mu_{\frak Q}([q],\la_1)+\mu_{\frak Q}([q],\la_2)$.
%%%
\subsubsection*{\q The global sections of $\O_{\frak P}(1,1)$}
%%%
All the assertions about the weights of one parameter subgroups
will rely on a good understanding of the global sections
of the line bundle $\O_{\frak P}(1,1)$ over ${\frak P}$.
Therefore, we will now give
an explicit description of them.
Let $[q\colon V\otimes \O_X\lra \E]$ be a point in ${\frak Q}$.
The points of $\widehat{\frak P}$ in the fibre over $[q]$ 
can then be written
as classes $[\eps, f]$, $\eps\in \C$, $f\in H:= 
\Hom(H^0(\E(\nu_0)),H^0(\E(\nu_0+m))^{\oplus u})$, 
and $[\eps, f]=[z\eps,zf]$, $z\in\C^*$.
Now, fix bases $v_1,...,v_p$ of $V$, $n_1,...,n_\nu$ of $N$,
and $m_1,...,m_\mu$ of $H^0(\O_X(m))$.
Using the lexicographic order, we get ordered bases for
$V\otimes N$ and $V\otimes N\otimes H^0(\O_X(m))$.
Set $p^\p:=P(\nu_0)$ and $p^{\p\p}:=P(\nu_0+m)$.
Let ${\frak I}$ be the set of all $p^\p$-tuples
of elements of the form $v_{\iota}\otimes n_{\kappa}$
whose images in $H^0(\E(\nu_0))$ form a basis.
Likewise ${\frak J}$ is defined as the set of
$p^{\p\p}$-tuples of the form $v_{\iota}\otimes n_{\kappa}\otimes m_\la$
which induce a basis for $H^0(\E(\nu_0+m))$.
Observe that specifying an element in $J\in {\frak J}$
is the same as specifying an element $S_J$ in 
$\bigwedge^{p^{\p\p}} V\otimes N\otimes H^0(\O_X(m))$ which,
viewed as a global section of ${\frak L}_{\frak Q}$,
does not vanish in $[q]$.
Pick elements $I\in {\frak I}$ and $J\in {\frak J}$,
and let $u_1,....,u_{p^\p}$ and $w_1,....,w_{p^{\p\p}}$
be the induced bases of $H^0(\E(\nu_0))$ and $H^0(\E(\nu_0+m))$,
respectively. We, thus, obtain a basis
$w_1^1,...,w^1_{p^{\p\p}}, ...,w_1^u,...,w^u_{p^{\p\p}}$ 
for $H^0(\E(\nu_0+m))^{\oplus u}$.
The $f_{ij}^k:= u_i^\vee \otimes w_j^k$, $i=1,...,p^\p$, $j=1,...,p^{\p\p}$,
$k=1,...,u$,
form a basis for $H$.
As one knows from linear algebra,
${w^{k\vee}_j}$ can be identified --- up to a sign --- with
$(w^k_1\wedge\cdots\wedge \hat{w}^k_j\wedge\cdots\wedge w^k_{p^\p})
/(w^k_1\wedge \cdots\wedge w^k_{p^\p})$.
Thus, $u_i\otimes w_j^{k\vee}$ 
defines a rational section $\sigma^{ijk}_{IJ}$
of 
$V\otimes N\otimes
\bigwedge^{p^{\p\p}-1} (V\otimes N\otimes M)\otimes {\frak L}^\vee_{\frak Q}$.
Therefore, $\widetilde{\Theta}^{ijk}_{IJ}:= \sigma^{ijk}_{IJ}\otimes S_J$
is a global section of 
$V\otimes N\otimes
\bigwedge^{p^{\p\p}-1} (V\otimes N\otimes M)\otimes\O_{\frak Q}$.
Denote the induced section of $\O_{\widehat{\frak P}}(1,1)$
by $\Theta_{IJ}^{ijk}$.
From the construction, it is clear that the $\Theta_{IJ}^{ijk}$
generate $\O_{\widehat{\frak P}}(1,1)$ everywhere and that they
are eigenvectors for the action of the maximal torus
defined by the basis $v_1,...,v_p$.
%%%
%\begin{Rem}
%The sets ${\frak I}$ and ${\frak J}$ inherit a lexicographic order.
%Therefore, we can always choose the minimal elements w.r.t.\
%that order.
%\end{Rem}
%%%
Let's return to the original setting.
Let $(\E,\eps,\phi,\psi)$ be a framed Hitchin pair.
Call a subsheaf $\F$ of $\E$ \it $\phi$-superinvariant\rm ,
if $\F\subset\ker\psi$ and the induced homomorphism
$\E/\F\lra (\E/\F)\otimes\G$ is also zero.
As an immediate consequence
of the previous discussion, we note
%%%
\begin{Lem}
\label{R1}
Fix a basis $v_1,...,v_p$ for $V$. Then one obtains the following values
for the action of $\la^{(i)}$ on the fibre of ${\frak N}_{\frak P}$
over $\lim_{z\rightarrow 0}{\bf p}\cdot\la^{(i)}(z)$
with ${\bf p}=([q\colon V\otimes\O_X\lra\E],[\eps,\phi])$, $i=1,...,p$:
{\rm i)} $-p$ if $\E_{\langle\, v_1,...,v_i\,\rangle}$ is not
$\phi$-invariant;
{\rm ii)} $p$ if $\eps=0$ and $\E_{\langle\, v_1,...,v_i\,\rangle}$ is
      $\phi$-superinvariant; and
{\rm iii)} $0$ in all the other cases.
\end{Lem}
%%%
The above description of the weights and a standard argument in
Simp\-son-type constructions --- \cite{HL} being closest to our
situation --- lead to the following conclusion:
%%%
\begin{Prop}
\label{R2}
Let $r=([q\colon V\otimes\O_X\lra\E],[\eps,\phi],[\psi])$ 
be in the closure of ${\frak R}_0$.
Then the following assertions are equivalent:
\begin{enumerate}
\item For any basis $v_1,...,v_p$ of $V$ and any one parameter
      subgroup $\la^{(i)}$, such that the sheaf $\E_{\langle\, v_1,...,v_i\,\rangle}$
      is $\phi$-invariant but not $\phi$-superinvariant when $\eps=0$,
      $\mu(r,\la^{(i)})(\ge) 0.$
\item For all subspaces $U\subset V$ such that  $\E_U$
      is $\phi$-invariant but not $\phi$-super\-in\-variant when $\eps=0$,
      \begin{eqnarray*}
      {\dim U\over\rk \E_U} - {\sigma(0)\over\rk \E_U} &(\le)&
      {p\over r} -{\sigma(0)\over r}\\
      {\dim U\over\rk \E_U}
      &(\le)&
      {p\over r} -{\sigma(0)\over r},\qquad\hbox{\rm if}\q
      \E_U\subset\ker\phi.
      \end{eqnarray*}
\end{enumerate}
\end{Prop}
Now, we proceed to the proof of Theorem~\ref{SemStab}.
Let's start with the following obvious 
%%%
\begin{Cor}
\label{curve}
Suppose $X$ is a curve and $r$ lies in $\overline{{\frak R}_0}$.
Then $H^0(q)$ must be an isomorphism and $\E$ must be torsion free.
\end{Cor}
%%%
\begin{Rem}
In the higher dimensional cases, one cannot adapt the proof for the 
moduli spaces of semistable sheaves.
First, one can copy the proof of Proposition~4.4.2 in \cite{HL2} 
to get a homomorphism $\kappa\colon \E\lra\E^\p$.
One can even equip $\E^\p$ with the structure of a 
framed Hitchin pair such that $\kappa$ becomes 
a homomorphism between framed Hitchin pairs.
But unfortunately, it is not clear whether $\ker(V\lra H^0(\E^\p))$
generates a $\phi$-invariant (torsion) subsheaf of $\E$.  
\end{Rem}
%%%
We will apply the Criterion~iii) of Proposition~\ref{Sect}.
First, assume that $r$ is a (semi)stable point.
Then, by Proposition~\ref{R2}, we only have to show that $\phi$ can't be
nilpotent when $\eps=0$. If $\eps=0$ and $\phi$ is nilpotent,
then there is a filtration
$$
0\subset \F_1\subset \cdots\subset \F_s\subset\F
$$
such that $\phi(\F_i)\subset \F_{i-1}\otimes\G$ and the $\F_i$
are globally generated, $i=1,...,s$ (Assumption~\ref{Ass1}).
Choose a basis $v_1,...,v_p$ such that there are indices
$\iota_1<\cdots<\iota_s$ with $\langle\, v_1,...,v_{\iota_i}\,\rangle
=H^0(\F_i)$, $i=1,...,s$. Let $\la$ be given w.r.t.\ that basis
by $\sum_{i=1}^s\gamma^{(\iota_i)}$. Then,
semistability and Assumption~\ref{Ass1} yield
$$
-p+2\sum_{i=1}^s \Bigl(p h^0(\F_i(\nu_0+m))-h^0(\F_i)P(\nu_0+m)+
{a\over 2}(p-h^0(\F_s))\Bigr)\ge 0.
$$
Plugging in our formula for $a/2$, viewing everything as polynomials
in $\nu_0$ and taking leading coefficients gives
$$
-p\ +\ \sigma(0)\ +\ 2\sum_{i=1}^s \Bigl(p\rk\F_i-h^0(\F_i)r+
\sigma(0)(r-\rk\F_i)\Bigr)\q \ge\q 0.
$$
If we have chosen $n_1$ in~\ref{Ass1} big enough, then this is not
a possibility. Indeed, $-P+\sigma$ --- as a function of $n$ ---
is a polynomial of degree $\dim X$ with negative leading coefficient 
whereas the polynomial
corresponding to the sum has at most degree $\dim X-1$.
Since, by Corollary~\ref{bound3}, there can occur only 
finitely many such polynomials, we are done.
\par
To see the converse, let $(\E,\eps,\phi,\psi)$ satisfy Condition~iii)
of Prop.~\ref{Sect}. Then, the second condition in
Proposition~\ref{R2} is satisfied.
Let $v_1,...,v_p$ be any basis and $\la$ be a one parameter subgroup
given by, say, $\sum_{i=1}^{p-1} \alpha_i\gamma^{(i)}$, $\alpha_i\in
\Z[(1/p)]_{\ge 0}$.
Note that, if $\phi$ is not nilpotent, there is a non-zero
global section of $\O(2,1,a)$ on which every one
parameter subgroup $\la$ acts with weight $\le -2\mu_{\frak Q}([q],\la)
-a\mu_{\bf R}([\psi],\la)$.
Therefore, if $\E_{\langle\, v_1,...,v_i\,\rangle}$
is $\phi$-invariant
for every $i$ with $\alpha_i\neq 0$,
there is nothing
to show.
Otherwise let $i_1,...,i_s$ be the indices belonging to
non-invariant subsheaves.
For each $j=1,...,s$, let $\E_{i_j}$ be the saturation of $\E_{V_{i_j}}$.
We then find indices $\iota_1,...,\iota_t$ among $i_1,...,i_s$
such that $\E_{V_{i_j}}\subset \E_{\iota_\kappa}$ if and only if
$i_j\le \iota_\kappa$, $j=1,...,s$, $\kappa=1,...,t$.
For $\iota_{\kappa-1}<i_j\le \iota_{\kappa}$, the 
induced homomorphism $\E_{V_{i_j}}\lra\E/\E_{\iota_{\kappa}}\otimes\G$
will be non-zero. Therefore, we find a section $\Theta_\kappa$
in $H^0({\frak P}, \O_{\frak P}(1,1))$ 
among the $\Theta_{IJ}^{ijk}$ such that $\la^{(i_j)}$ acts on
$\Theta_\kappa$ with weight $-p-\mu_{\frak Q}([q],\la^{(i_j)})$ 
for every $j$ with $\iota_{\kappa-1}<i_j\le \iota_{\kappa}$ and any other
one parameter group $\la^{(i)}$ with weight
$\le -\mu_{\frak Q}([q],\la^{(i)})$. 
Let $\Theta:=\Theta_1\otimes\cdots\otimes \Theta_t$ be the corresponding
section of $\O_{\frak P}(t,t)$. Then $\la$ acts on $\Theta$
with weight $\le -p\sum_{j=1}^s \alpha_{i_j}- t\sum_{i=1}^p\alpha_i
\mu_{\frak Q}([q],\la^{(i)})$.
Therefore, the assertion $\mu(r,\la)(\ge) 0$ can be reduced
to
$$
p \q +\q 2t\bigl(\mu_{\frak Q}([q], \la^{(i_j)})+
{a\over 2}\mu_{\bf R}([\psi],\la^{(i_j)})\bigr)
\q >\q0,\qquad j=1,...,s.
$$
Now, only those $i_j$ matter for which
$\mu_{\frak Q}([q], \la^{(i_j)})+
(a/ 2)\mu_{\bf R}([\psi],\la^{(i_j)})$ is negative.
But then, $\E_{V_{i_j}}$ can be assumed to be globally generated and without
higher cohomology, and replacing this sheaf with its saturation,
we are left to show
$$
p\ -\ \sigma(0)\ +\ 2t\Bigl(\rk\E_{\iota_\tau} (p-\sigma(0))- r h^0(\E_{\iota_\tau})\Bigr)
\q >\q 0,\qquad \tau=1,...,t.
$$
We view this again as an inequality between polynomials in $n$.
By Remark~\ref{PolVersion},
the second term is then $\ge - 2t(r-1)[rC^\p n^{\dim X-1}+\sigma(n)]$,
i.e., bounded from below by 
a polynomial of degree $\dim X-1$, and $P(n)-\sigma(n)$ is a polynomial
of degree $\dim X$ with positive leading coefficient, whence the claim.
\par
For the assertion about the polystable points,
let the one parameter subgroup $\la$ be
given by $\sum_{i=1}^{p-1} \alpha_i\gamma^{(i)}$, $\alpha_i\in
\Z[(1/p)]_{\ge 0}$, w.r.t.~the basis $v_1,...,v_p$ of $V$.
Observe that the above proof shows that
$\mu(\la, r)=0$ can only occur if for each $\alpha_i\neq 0$,
$\langle\, v_1,...,v_i\,\rangle=:V_i=H^0(\E_{V_i})$,
$\E_{V_i}$ is $\phi$-invariant
and destabilizes the framed Hitchin pair
$(\E,\eps,\phi,\psi)$.
%%%
So, $r$ will be a fixed point for the action of any such $\la$ if and
only if
$(\E,\eps,\phi,\psi)$
is a $\sigma$-polystable framed Hitchin pair.
$\square$
%%%
\subsection*{\q\sc Some variations and examples}
\label{Var}
%%%
In this section, let $X$ be a curve, and $L=\G$ a line bundle.
The type will be written as $(d,r,L,\H)$.
%%%
\subsubsection*{\q Framed Hitchin pairs with fixed determinant}
%%%
Fix a line bundle $M$ of degree $d$. A \it framed Hitchin pair
of type $(M,r,L,\H)$ \rm  is a framed Hitchin pair $(E,\eps,\phi,\psi)$
of type $(d,r,L,\H)$ with $\det(E)\cong M$.
The equivalence relation is the same as for framed Hitchin pairs
of type $(d,r,L,\H)$.
Observe that the sheaf ${\frak E}_{\frak R}$ on ${\frak R}$
provides us with a morphism $d({\frak E}_{\frak R})\colon {\frak R}
\lra \Pic X$. Since this morphism is $\SL(V)$-invariant,
it descends to a map
${\frak d}\colon {\cal {\cal F}{\cal H}}^{\sigma-ss}_{d/r/L/\H}\lra
\Pic X$. Define ${\cal {\cal F}{\cal H}}^{\sigma-ss}_{M/r/L/\H}$ as the scheme
theoretic fibre of ${\frak d}$ over $[M]$.
%%%
\subsubsection*{\q 
The condition ${\mathop{\rm Im}}(\phi)\subset\ker(\psi)\otimes L$}
%%%
Stupariu looks at framed Hitchin pairs which satisfy
$(\psi\otimes \id_L)\circ \phi=0$, i.e.,
$$
{\mathop{\rm Im}}\phi \qquad\subset\qquad (\ker \psi)\otimes L.
$$
Let us call these objects
\it framed Hitchin pairs of type $(d,r,L,\H,*)$\rm .
The definitions of \it equivalence\rm , \it $\sigma$-(semi)stability\rm ,
and so on carry over.
Note that the above condition forces $\ker\psi$ to be $\phi$-invariant.
We can also construct moduli spaces for those objects:
Consider, on ${\frak R}\times X$,
$$
\pi_X^*\H^\vee
\stackrel{\psi_{{\frak R}}^\vee}{\lra}
{\frak E}_{{\frak R}}^\vee
\stackrel{
\phi_{{\frak R}}^\vee\otimes {\id}_{\pi_{{\frak R}}^*{\frak N}_{{\frak R}}}\otimes {\id}_{\pi_X^*L}}{\lra}
{\frak E}_{{\frak R}}^\vee\otimes {\pi_{{\frak R}}^*{\frak N}_{{\frak R}}}\otimes\pi_X^*L.
$$
Define ${\frak R}^*$ as the closed subscheme of ${\frak R}$
whose closed points are the points $r\in {\frak R}$
such that the above homomorphism becomes zero when restricted
to $\{r\}\times X$.
One can now go on as before.
We denote the resulting moduli spaces
by ${\cal {\cal F}{\cal H}}_{d/r/L/\H/*}^{\sigma-(s)s}$.
Moreover, we can fix a line bundle $M$ and define the moduli
spaces  ${\cal {\cal F}{\cal H}}_{M/r/L/\H/*}^{\sigma-(s)s}$.
%%%
\subsubsection*{\q Some observations}
%%%
Let $(E,\eps,\phi,\psi)$ be
of type $(d,2,\allowbreak L,\O_X(m_0),*)$, $m_0\in\Z$.
The image of $\psi$ is of the form $\O_X(m_0)(-D)$ for
some effective divisor $D$, and we get an extension
$$
(e)\colon 0\lra \det(E)(D)(-m_0)\lra E\lra \O_X(m_0)(-D)\lra 0.
$$
%%%
\begin{Lem}
\label{Obs}
{\rm i)} If $\sigma>2m_0-d$, then there are no $\sigma$-semistable
framed Hitchin pairs of type $(d,2,L,\O_X(m_0),*)$.
\par
{\rm ii)} If $(E,\eps,\phi,\psi)$  is a $\sigma$-(semi)stable
Hitchin pair of type $(d,2,L,\allowbreak\O_X(m_0),\allowbreak *)$,
then there exists a $\sigma^\p \in {\Bbb Q}_{> 0}$ such that
$(E,\psi)$ is a $\sigma^\p$-(semi)stable framed module in the sense
of \cite{HL}.
\end{Lem}
%%%
\it Proof\rm.
i) Indeed, if $E$ is given as an extension $(e)$ as before,
then $\sigma$-semi\-stability implies $d-m_0\le (d-\sigma)/2$,
because $\ker\psi$ is $\phi$-invariant.
We also see that we must choose $m_0> d/2$.
\par
ii) Again, let $(E,\psi)$ be given by the extension $(e)$.
We know $d+\deg(D)- m_0\le (d-\sigma)/2$.
Set
$
\sigma^\p:=-d-2\deg(D)+2m_0.
$
Then $d+\deg(D) - m_0 = (d-\sigma^\p)/2$,
and for any line subbundle $F\neq \ker\psi$ of $E$ we have
$\deg(F)\le -\deg(D)+m_0 = (d+\sigma^\p)/2$.
If the extension is non-split, then we can choose $\sigma^\p$
slightly smaller, so that $(E,\psi)$ is even $\sigma^\p$-stable.
$\square$\par\bigskip
%%%
\begin{Rem}
The analogous statement of~\ref{Obs}~i) for framed Hitchin pairs
of type $(d,2,L,\H)$ is false. Indeed, consider, e.g.,
$X=\P_1$, $L=\H=\O_{\P_1}$, $E=\O_{\P_1}^{\oplus 2}$, $\eps:=1$, 
let $\psi$ be the projection onto the first factor
and $\phi$ be given by the matrix
$$
\left(
\begin{array}{cc}
 0 & 1\\
 0 & 0\\
\end{array}
\right).
$$
Then $(E,\eps,\phi,\psi)$ is $\sigma$-stable
for any $\sigma>0$.
We will come back to this in the section about boundedness in the next chapter.
\end{Rem}
%%%
\begin{Ex}
\label{Ex1}
We look at the situation
$X=\P_1$, $d=0$, $r=2$, $L=\H=\O_{\P_1}$.
Then, for a $\sigma$-semistable framed Hitchin pair $(E,\eps,\phi,\psi)$,
we will have $E\cong\O_{\P_1}\oplus \O_{\P_1}$, and $\ker\psi$ can't be
$\phi$-invariant.
Thus, $(E,\eps,\phi,\psi)$ is $\sigma$-semistable if and only
if $\ker\psi$ is not $\phi$-invariant and either $\eps\neq 0$ or
$\phi^2\neq 0$.
In particular, $(E,\eps,\phi,\psi)$ is then stable for any $\sigma>0$,
and $(E,\eps,\phi)$ is a semistable Hitchin pair.
Denote by ${\cal {\cal F}{\cal H}}$ the moduli space, and by ${\cal H}$ the
moduli space of semistable Hitchin pairs. There
is a natural map $\pi\colon {\cal {\cal F}{\cal H}}\lra {\cal H}$.
The space ${\cal H}$ is isomorphic to $\P_2$ with
coordinates, say, $[l_0,l_1,l_2]$.
Here, the point $[l_0,l_1,l_2]$ represents the class of the Hitchin pair
$(\O_{\P_1}\oplus \O_{\P_1}, \eps, \phi)$ with
$$
\eps=l_0,\q \phi= \left(\begin{array}{cc}
                          l_1 & 0\\
                           0  & l_2\\
                         \end{array}\right).
$$
It is easy to describe the fibres of $\pi$:
The preimage of $[l_0,l_1,l_2]$ with $l_1\neq l_2$ consists
just of the class of $(\O_{\P_1}\oplus \O_{\P_1}, \eps, \phi, \psi)$
with
$$
\eps=l_0,\q \phi= \left(\begin{array}{cc}
                          l_1 & 0\\
                           0  & l_2\\
                         \end{array}\right),\q
\psi= (1,1);
$$
and the preimage of $[l_0,l,l]$
of the class of $(\O_{\P_1}\oplus \O_{\P_1}, \eps, \phi, \psi)$
with
$$
\eps=l_0,\q \phi= \left(\begin{array}{cc}
                          l & 1\\
                           0  & l\\
                         \end{array}\right),\q
\psi= (1,0).
$$
Hence,  $\pi$ is an isomorphism.
It is possible to give explicit coordinates for ${\cal {\cal F}{\cal H}}$.
For this, we write
${\cal {\cal F}{\cal H}}=\P(M_1^\vee\oplus M_2^\vee)\catqot(\C^*\times \SL_2(\C))$,
where $M_1=\C\oplus \C^{(2,2)}$ and $M_2={\C^2}^\vee$, and $\C^*$
acts with weights $-1$ and $2$, so that the induced
polarization on $\P(M_1^\vee)\times\P(M_2^\vee)$ is $\O(2,1)$.
We choose coordinates $(l_0,l_{1,1},l_{1,2},l_{2,1},l_{2,2},
s_1,s_2)$. Observe, that a $(2\times 2)$-matrix is nilpotent
if and only its determinant and its trace are both $=0$.
Moreover, ``$\ker\psi$ is $\phi$-invariant'' can be expressed
as
\begin{eqnarray*}
D &:=& \det\left(\begin{array}{cc} l_{1,1} s_2 - l_{2,1}s_1 & s_2\\
                                   l_{2,1} s_2 - l_{2,2}s_1 & s_1
              \end{array}\right)\\
  &=& s_1s_2(l_{1,1}+l_{2,2}) - s_2^2 l_{2,1}- s_1^2 l_{1,2} =0.
\end{eqnarray*}
Thus, the $\SL_2(\C)$-nullforms are the common zeroes
of the polynomials $H_0:=l_0$,
$H_1:= l_{1,1}+l_{2,2}$,
$H_2:= l_{1,1}l_{2,2}-l_{1,2}l_{2,1}$, and $H_3:=D$.
We have to find those homogeneous polynomials
in $H_0,...,H_3$ which are $\C^*$-invariant. The weights of
$H_0$, $H_1$, $H_2$, and $H_3$ w.r.t.\ the $\C^*$-action are $1$, $1$, $2$
and $-2$. Hence, the coordinates are given by
$h_0=H_0^2H_4$, $h_1=H_1^2H_4$, and $h_3=H_3H_4$.
\end{Ex}
%%%
\section{Oriented framed Hitchin pairs}
%%%
We will first discuss the notion of oriented framed Hitchin pairs
and introduce the --- parameter independent --- semistability concept
for them.
Then we will proceed to construct the moduli spaces of semistable
oriented framed Hitchin pairs over curves. There are
some intricate technical points such as the behaviour of
$\sigma$-semistability when $\sigma$ becomes large.
%%%
\subsection*{\q\sc Oriented (symmetric) framed Hitchin pairs}
%%%
An \it oriented fra\-med Hitchin pair of type $(P,\G,\H,\n)$ \rm is a
quintuple $(\E,\eps,\delta,\phi,\psi)$ where $\E$, $\eps$, $\phi$,
and $\psi$ have the same meaning as before, only that $\psi=0$ is now
allowed, and $\delta\colon\det\E\lra{\cal N}[\E]$ is a homomorphism.
%%%
An \it isomorphism between oriented framed Hitchin pairs
$(\E,\eps,\delta,\phi,\psi)$  and
$(\E^\p,\eps^\p,\delta^\p,\phi^\p,\psi^\p)$ of type $(P,\G,\H,\n)$ \rm
is an isomorphism $\rho\colon \E\lra\E^\p$ such that there exist
numbers $w,z\in\C^*$ such that
$$
\eps^\p=z\eps,\ \delta^\p=w^r \delta\circ(\det\rho)^{-1},
\ \phi^\p=z\bigl((\rho\otimes{\id}_{\G})\circ\phi\circ\rho^{-1}\bigr),
\ \hbox{\& } \psi^\p=w\psi\circ\rho^{-1}.
$$
An isomorphism $\rho$ as above will be called a \it proper
isomorphism\rm,
if $w=1$. Note that both notions of isomorphism yield the
same equivalence relation on the set of all oriented framed Hitchin
pairs of type $(P,\G,\H,\n)$.
We will call an oriented framed Hitchin pair
$(\E,\eps,\delta,\phi,\psi)$ of type $(P,\G,\H,\n)$ \it symmetric\rm ,
if $\phi$ is symmetric. As before, any symmetric oriented framed Hitchin pair
of type $(P,\G,\H,\n)$ defines an element
$\widehat{\chi}(\E,\eps,\delta,\phi,\psi)\in\widehat{\P}_\G$.
This depends only on its equivalence class and is called the
\it characteristic polynomial of $(\E,\eps,\delta,\phi,\psi)$\rm .
In general. we can assign to every oriented framed Hitchin pair
$(\E,\eps,\delta,\phi,\psi)$ its \it characteristic vector \rm in 
$\widehat{\Xi}_\G$.
%%%
\begin{Rem}
\label{Aut}
i) The automorphisms of the oriented framed Hitchin pairs living in the
universal family which will be constructed below coming from actions
of the stabilizers are only automorphisms in the weaker sense.
Thus, one should adopt this notion of isomorphism in order to
avoid confusion when trying to make the universal family descend.
The above notion of isomorphism is also the one which extends to
families. This situation is unlike the situation for framed Hitchin
pairs! Note that $\la\cdot\id$, $\la\in\C^*$, is always an automorphism.
Moreover, if $\rho$ is an isomorphism with constant $w$, then
$w^{-1}\rho$ will be a proper isomorphism.
\par
ii) If $X$ is a curve and $L=\G$ is a line bundle, we can also extend
the definitions of Chapter~\ref{Var}, i.e., we can define the notion
of an \it oriented framed Hitchin pair of type $(d,r,L,\H,\n,*)$\rm,
\it of type $(M,r,L,\H)$ \rm and $(M,r,L,\H,*)$. In the latter two cases
we mean $\det(E)\cong M$, and $\delta\colon\det(E)\lra M$.
\end{Rem}
%%%
A \it family of oriented framed Hitchin pairs of type $(P,\G,\H,\n)$
para\-metri\-zed by the noetherian scheme $S$ \rm is defined to be a
seventuple
$({\frak
E}_{S},\eps_S,\delta_S,\allowbreak \phi_S,\widehat{\psi}_S,{\frak M}_S,
{\frak N}_S)$ consisting of line bundles ${\frak M}_S$ and ${\frak N}_S$
on $S$, an $S$-flat family of torsion free coherent sheaves ${\frak
E}_S$ with Hilbert polynomial $P$ on $S\times X$, a section
$\eps_S\in H^0({\frak N}_S)$, a homomorphism $\delta_S\colon \det({\frak
E}_S)\lra {\cal N}[{\frak E}_{S}]\otimes \pi_S^*{\frak M}_S$,
a twisted endomorphism $\phi_S\colon {\frak E}_S\lra
{\frak E}_S\otimes\pi_S^*{\frak N}_S\otimes\pi_X^*\G$, and a
homomorphism $\widehat{\psi}\colon S^r {\frak E}_S \lra \pi_X^*
S^r\H\otimes \pi_S^*{\frak M}_S$ which is outside the closed
subscheme $S_0:=\{\, s\in S\,|\,\widehat{\psi}_{S|\{s\}\times X}=0\,\}$
a symmetric power, i.e., there is a line bundle ${\frak M}^\p$
on $S\setminus S_0$ with ${\frak M}^{\p\otimes r}={\frak
M}_{S|S\setminus S_0}$, and
$\widehat{\psi}_{S|S\setminus S_0}$
is the symmetric
power of a homomorphism $\psi^\p\colon {\frak E}_{S|(S\setminus
S_0)\times X} \lra
\pi_X^*\H\otimes\pi_{S\setminus S_0}^*{\frak
M}^\p$.
We leave it to the reader to define \it equivalence of families\rm.
%%%
\subsection*{\q\sc Semistability}
The definition of semistability will be made in analogy to the
definition of semistability for oriented pairs in \cite{OST}.
We need a preparatorial result.
%%%
\begin{Lem}
\label{HN1}
Let $(\E,\eps,\phi,\psi)$ be a framed Hitchin pair of type $(P,\G,\H)$,
possibly with $\psi=0$.
Suppose there are non-trivial $\phi$-invariant subsheaves in $\ker\psi$.
Then there exists a uniquely determined non-trivial $\phi$-invariant
subsheaf
$\K_{\max}\subset \ker\psi$, such that for all other $\phi$-invariant
subsheaves $\K\subset\ker\psi$, one has $(P(\K)/\rk\K)\le (P(\K_{\max})/
\rk\K_{\max})$, and, if equality occurs, $\K\subset\K_{\max}$.
\end{Lem}
%%%
\it Proof\rm.
Indeed, since, by assumption, the set of non-trivial $\phi$-invariant
subsheaves
in $\ker\psi$ is not empty and the sum and the intersection of two
$\phi$-invariant subsheaves in $\ker\psi$ is again a $\phi$-invariant
subsheaf in $\ker\psi$, to get the result, one merely needs to copy the
proof of Lemma~1.3.5 in~\cite{HL2}.
$\square$\par\bigskip
%%%
Set
$
\sigma_{\E,\phi,\psi}:= P(\E) -
(\rk\E/\rk\K_{\max})P(\K_{\max}).
$
Let $(\E,\eps,\delta,\phi,\psi)$ be an oriented framed Hitchin pair of
type $(P,\G,\H,\n)$. We call it \it semistable\rm,
if and only if either
there are no $\phi$-invariant subsheaves in
$\ker\psi$, or $\delta$ is an isomorphism  and there are
$\phi$-invariant subsheaves in
$\ker\psi$, $\sigma_{\E,\phi,\psi}\ge 0$ and for all non-trivial
$\phi$-invariant subsheaves $\F\subset\E$
$$
{P_\F\over \rk\F}-{\sigma_{\E,\phi,\psi}\over\rk\F}\q\le\q
{P_\E\over \rk\E}-{\sigma_{\E,\phi,\psi}\over\rk\E}.
$$
And we call $(\E,\eps,\delta,\phi,\psi)$ \it stable\rm, if and only if
either there are no $\phi$-invariant subsheaves in
$\ker\psi$, or $\delta$ is an isomorphism  and there are
$\phi$-invariant subsheaves in
$\ker\psi$, $\sigma_{\E,\phi,\psi}> 0$, and one of the following two
possibilities holds:
\begin{enumerate}
\item For all non-trivial
$\phi$-invariant proper subsheaves $\F\subset\E$
$$
{P_\F\over \rk\F}-{\sigma_{\E,\phi,\psi}\over\rk\F}\q<\q
{P_\E\over \rk\E}-{\sigma_{\E,\phi,\psi}\over\rk\E}.
$$
\item $\psi\neq 0$, $(\E,\eps,\phi,\psi)$ splits
      as $(\K_{\max},\eps,\phi_{|\K_{\max}},0)\oplus
       (\E^\p,\eps,\phi_{|\E^\p},\psi)$,
      the triple $(\K_{\max},\allowbreak \eps,\phi_{|\K_{\max}})$ is a stable Hitchin pair,
      and $(\E^\p,\eps,\phi_{|\E^\p},\psi)$ is a
      $\sigma_{\E,\phi,\psi}$-stable framed Hitchin pair, such that
$P(\K_{\max})/\rk\K_{\max}=(P(\E^\p)-\sigma_{\E,\phi,\psi})/\rk\E^\p$.
\end{enumerate}
For our constructions, we have to restate the semistability concept
in terms of the semistability concepts for framed Hitchin pairs of
Chapter~1.
%%%
\begin{Lem}
{\rm i)} An oriented framed Hitchin pair
$(\E,\eps,\delta,\phi,\psi)$ of type $(P,\allowbreak\G,\allowbreak\H,\n)$ is
semistable, if and only if it satisfies one of the following
three conditions:
\begin{enumerate}
\item There are no $\phi$-invariant subsheaves in the kernel of $\psi$.
\item $\delta$ is an isomorphism, and $(\E,\eps,\phi)$ is a semistable Hitchin 
      pair of type
      $(P,\G)$. 
\item $\psi\neq 0$, $\delta$ is an isomorphism, and there is a polynomial
      $\sigma\in{\Q}[t]$
      of degree less than $\dim X$ with positive leading coefficient
      s.~th.\
      $(\E,\eps,\phi,\psi)$ is a $\sigma$-semistable framed Hitchin
      pair of
      type $(P,\G,\H)$.
\end{enumerate}
{\rm ii)}
$(\E,\eps,\delta,\phi,\psi)$ is
stable, if and only if
satisfies one of the conditions listed below.
%%%
\begin{enumerate}
   \item There are no $\phi$-invariant subsheaves in $\ker \psi$.
   \item $\delta$ is an isomorphism, and 
         $(\E,\eps,\phi)$ is a stable Hitchin pair of type
         $(P,\G)$.
   \item $\psi\neq 0$, $\delta$ is an isomorphism, 
         and there is a polynomial $\sigma\in{\Q}[t]$
         of degree less than $\dim X$ with positive leading coefficient
         s.~th.\
         $(\E,\eps,\phi,\psi)$ is a $\sigma$-stable framed Hitchin pair
         of type $(P,\G,\H)$.
   \item $\psi\neq 0$, $\delta$ is an isomorphism, 
         and there is a polynomial $\sigma\in{\Q}[t]$
         of degree less than $\dim X$ with positive leading coefficient
         s.~th.\
         $(\E,\eps,\phi,\psi)$ is a $\sigma$-polystable framed Hitchin
         pair
         of type $(P,\G,\H)$ of the form $(\E^\p,\eps,\phi^\p,0)\oplus
         (\E^{\p\p},\eps,\phi^{\p\p},\psi)$.
\end{enumerate}
\end{Lem}
%%%
\begin{Rem}
\label{Obs100}
    In the case the base $X$ is a curve, $L=\G$ a line bundle, and $r=2$,
    then,
    by
    Lemma~\ref{Obs}, for
    any (semi)stable oriented framed Hitchin pair
    $(E,\eps,\delta,\phi,\psi)$
    of type $(d,2,L,\H,\n,*)$, then either $(E,\eps,\phi)$ is a
    (semi)stable Hitchin pair of type $(d,2,L)$, or the triple
    $(E,\delta,\psi)$ is a
    (semi)stable oriented pair of type $(d,2,\n)$ in the sense of
    \cite{OST}.
\end{Rem}
%%%
Let $(\E,\eps,\delta,\phi,\psi)$ be a semistable oriented framed Hitchin
pair of type $(P,\G,\H,\n)$ which is not stable. This occurs if and
only if there is either a $\phi$-invariant subsheaf in $\ker\psi$
which destabilizes $(\E,\eps,\phi)$ as a Hitchin pair --- in which
case $(\E,\eps,\phi)$ must be a semistable Hitchin pair --- or
there are a $\phi$-invariant subsheaf $\K$ of $\ker\psi$ and a
$\phi$-invariant subsheaf $\F\not\subset\ker\psi$ with
$$
{P_\F \over\rk\F}-{\sigma_\K\over\rk\F}\q =\q
{P_\E \over\rk\E}-{\sigma_\K\over\rk\E},
$$
where $\sigma_\K:=P(\E)-(\rk\E/\rk\K)P(\K)$.
Note that in this case $(\E,\eps,\phi,\psi)$ is $\sigma_\K$-semistable
but not $\sigma$-semistable for every polynomial $\sigma\neq\sigma_\K$.
Let $0\subset \E_m\subset\cdots\subset\E_1\subset \E$  be either the
Jordan-H\"older filtration of the semistable
Hitchin pair $(\E,\eps,\phi)$ or the Jordan-H\"older filtration of the
$\sigma_\K$-semistable framed Hitchin pair $(\E,\eps,\phi,\psi)$.
Note that such a filtration induces a canonical isomorphism
between $\det(\bigoplus_{i=1}^{m+1}\E_{i-1}/\E_i)$ and $\det(\E)$.
Hence, we obtain an  associated graded object
$\gr(\E,\eps,\delta,\phi,\psi)$ --- well-defined up to equivalence ---,
and we call
$(\E,\eps,\delta,\phi,\psi)$ \it polystable\rm, if it is equivalent
to its associated graded object. Furthermore, two semistable oriented
framed Hitchin pairs are said to be \it S-equivalent\rm, if and only
if their associated graded objects are equivalent.
%%%
\begin{Rem}
\label{Obs200}
i) At this moment, one could get the impression that a stable
   oriented
   framed Hitchin pair $(\E,\eps,\delta,\phi,\psi)$ where $(\E,\eps,
   \phi,\psi)$ splits into two components $(\E^\p,\eps,\phi^\p,\allowbreak 0)$ and
   $(\E^{\p\p},\eps,\phi^{\p\p},\psi)$
   might not be stable at all, because there might be another semistable
   oriented framed Hitchin pair degenerating to it. This is, of course,
   not the case. Let $\sigma_0$ be the unique polynomial w.r.t.\ which
   $(\E,\eps,\phi,\psi)$ is semistable. Suppose $(\widetilde{\E},\eps,
   \widetilde{\psi},\widetilde{\psi})$ is a $\sigma_0$-semistable
   pair whose associated graded object is equivalent to
   $(\E,\eps,\phi,\psi)$. Then it is easy to see that either it is
   itself equivalent to $(\E,\eps,\phi,\psi)$ or it is $\sigma$-stable
   w.r.t.\ some polynomial $\sigma$ which is "close" to $\sigma_0$.
\par
ii) If $(\E,\eps,\delta,\phi,\psi)$ is properly semistable 
    and $(\E,\eps,\phi)$ is not a semi\-stable Hitchin pair, then
    $m$ in the above Jordan-H\"older filtration must be at least
    two.
\par
iii) One checks that the stable oriented framed Hitchin pairs are
    precisely the polystable oriented framed Hitchin pairs which have
    only finitely many proper automorphisms.
\end{Rem}
%%%
\subsection*{\q\sc Boundedness}
%\label{Sec:Bound}
In this section, we prove the boundedness of the set of isomorphy
classes of torsion free coherent sheaves occuring in semistable
oriented framed Hitchin pairs of type $(P,\G,\H,\n)$ and carefully
examine how the notion of $\sigma$-semistability behaves when
$\sigma$ becomes in a certain sense large.
Invoking Maruyama's boundedness result \cite{Ma} again, the
boundedness will follow from
%%%
\begin{Prop}
\label{bound10} Suppose $\G=\O_X(m)^{\oplus u}$
and $\O_X(m)$ is globally generated. Then,
there is a constant $C$ such that, for any semistable oriented framed
Hitchin pair $(\E,\eps,\delta,\phi,\allowbreak \psi)$ of type $(P,\G,\H,\n)$, the
condition 
$$\mu_{\max}(\E)\le C$$ 
holds true.
\end{Prop}
%%%
\it Proof\rm. For $\psi=0$, this is Proposition~2.2.2 in \cite{OST}. Thus, we can
assume $\psi\neq 0$. Any subsheaf $\F$ of $\E$ can be written as an
extension
$$
0\lra \ker\psi\cap\F\lra\F\lra \psi(\F)\lra 0.
$$
Since $\mu(\psi(\F))$ is bounded by $\mu_{\max}(\H)$, it suffices
to bound $\mu_{\max}(\ker\psi)$.
Recall the following
%%%
\begin{Lem}
\label{Hom0}
Given torsion free coherent sheaves $\F_1$ and $\F_2$
with $\mu_{\min}(\F_1)>\mu_{\max}(\F_2)$. Then there does not exist
any non-trivial homomorphism  from $\F_1$ to $\F_2$.
\end{Lem}
%%%
Set $\mu_0:=\max\{\,\mu,\mu_{\max}(\H)\,\}$. We will derive a
contradiction from the following assumption:
%%%
$$
\mu_{\max}(\ker\psi)\q >\q \mu_0 + r\deg\O_X(m).
$$
%%%
Suppose this was true and let
$0=\K_0\subset\K_1\subset\cdots\subset \K_\kappa=\ker\psi$ be the
slope Harder-Narasimhan filtration.
%%%
\begin{Claim} For $i=1,...,\kappa-1$, the following inequality is
satisfied:
$$
\mu(\K_i/\K_{i-1})\q \le\q \mu(\K_{i+1}/\K_i)+\deg\O_X(m),
$$
in particular,
$$
\mu_{\max}(\ker\psi)\q \le\q \mu_{\min}(\ker\psi)+(r-1)\deg\O_X(m).
$$
\end{Claim}
%%%
Therefore, $\mu(\ker\psi)>\mu$, so that $\ker\psi$ cannot be
$\phi$-invariant. On the other hand,
\begin{eqnarray*}
\mu_{\min}(\ker\psi)& > &
\mu_{\max}(\H)+\deg\O_X(m)\q=\q\mu_{\max}(\H\otimes\O_X(m)^{\oplus u})\\
&\ge& \mu_{\max}((\E/\ker\psi)\otimes \O_X(m)^{\oplus u}).
\end{eqnarray*}
%%%
In view of Lemma~\ref{Hom0}, this means that $\ker\psi$ must be
$\phi$-invariant.
\par
To see the claim, first note that
$\phi(\K_1)\subset\ker\psi\otimes\O_X(m)^{\oplus u}$, again by~\ref{Hom0}. 
By semistability, $\K_1$ cannot be $\phi$-invariant.
Hence, there is an index $i^\p>1$, s.~th.\ $\phi(\K_1)\subset
\K_{i^\p}\otimes \O_X(m)^{\oplus u}$ and
$\phi(\K_1)\not\subset
\K_{i^\p-1}\otimes \O_X(m)^{\oplus u}$. Thus, there is a
non-trivial homomorphism $f\colon \K_1\lra
(\K_{i^\p}/\K_{i^\p-1})\otimes\O_X(m)^{\oplus u}$, whence
$$
\mu(\K_1)\le \mu(\K_{i^\p}/\K_{i^\p-1})+\deg\O_X(m)
\le \mu(\K_{1}/\K_2)+\allowbreak \deg\O_X(m).
$$
Next, suppose the claim is true for $i=1,...,j$.
Since $\mu_{\min}(\K_{j+1})=\mu(\K_{j+1}/\K_j)\ge
\mu(\K_1)-j\deg\O_X(m)$,
again $\phi(\K_{j+1})\subset \ker\psi\otimes\O_X(m)^{\oplus u}$,
so that the same argumentation as before goes through, and we settle the
case $j+1$.
$\square$\par\bigskip
%%%
\begin{Cor}
\label{bound20}
Let $(\E,\eps,\phi,\psi)$ be a framed Hitchin pair of type
$(P,\G,\allowbreak\H)$, such that there is no $\phi$-invariant
subsheaf which is contained in $\ker\psi$, then $(\E,\eps,\phi,\psi)$
will be $\sigma$-stable for all polynomials $\sigma\in{\Q}[t]$
of degree $\dim X-1$ with sufficiently positive leading coefficient.
\end{Cor}
%%%
\it Proof\rm.
In the above proof, we have ruled out that one of the $\K_i$ be
$\phi$-invariant by the semistability condition, here, we do it by
assumption. Thus, the same conclusion as in Theorem~\ref{bound10}
--- with the same constant $C$ --- holds for framed Hitchin pairs
$(\E,\eps,\phi,\psi)$ with non-trivial framing and no $\phi$-invariant
subsheaves in $\ker\psi$.
So, any polynomial
$\sigma\in{\Q}[t]$
of degree $\dim X-1$ with leading coefficient $>r(r-1)C-(r-1)d$ will do
the trick.
$\square$\par\bigskip
%%%
We
also have the converse
%%%
\begin{Prop}
\label{bound30}
For all polynomials $\sigma$ of degree $\dim X-1$
whose leading coefficient is sufficiently large and all
$\sigma$-semistable framed Hitchin pairs $(\E,\eps,\phi,\psi)$
of type $(P,\G,\H,\n)$, there will be no $\phi$-invariant
subsheaf in $\ker\psi$.
\end{Prop}
%%%
\it Proof\rm.
Let $(\E,\eps,\phi,\psi)$ be a framed Hicthin pair, set
$\F_0:=\ker\psi$,
and for $i\ge 1$, $\F_i:=\ker\bigl(\F_{i-1}\lra
(\E/\F_{i-1})\otimes\O_X(m)^{\oplus u}\bigr)$.
This yields a decreasing chain of saturated submodules
$$
0\subset\cdots\subset\F_{i+1}\subset\F_i\subset\cdots\subset\F_0.
$$
By definition, a $\phi$-invariant subsheaf $\F\subset\ker\psi$ is
contained in all the $\F_i$. Therefore, such a subsheaf exists
if and only if one of the $\F_i$ is $\phi$-invariant.
But, by construction, the $\F_i$'s coming from framed Hitchin pairs
which are $\sigma$-semistable for some polynomial $\sigma$ form
bounded families.
This means, if $[d-$(leading coefficient of $\sigma$)$]/r$ is smaller
than every possible $\mu(\F_i)$, then a $\sigma$-semistable framed
Hitchin pair $(\E,\eps,\phi,\psi)$ of type $(P,\G,\H,\n)$ has no
non-trivial $\phi$-invariant subsheaves which are contained in
$\ker\psi$.
$\square$\par\bigskip
%%%
\subsection*{\q\sc Flips between the moduli spaces of framed Hitchin pairs}
%\label{Flips}
%%% 
We will now carry out a
discussion of this topic which makes the phenomenon completely
transparent, and which is independent of the existence of master spaces.
\par
Let ${\Q}[t]_{\dim X-1, +}$ be the set of all polynomials
of degree at most $\dim X-1$ which have positive leading coefficient.
This set is totally ordered by the lexicographic order.
Let ${\frak OFH}_{P/\G/\H/\n}^{ss}$ be the set of equivalence classes
of semistable oriented framed Hitchin pairs of type $(P,\G,\H,\n)$.
Given the equivalence class of a semistable framed oriented Hitchin pair
$(\E,\eps,\delta,\phi,\psi)$,
then a non-trivial $\phi$-invariant saturated subsheaf $\K\subset\ker\psi$
of $\E$
defines a polynomial
$\sigma_\K\in {\Q}[t]_{\dim X-1, +}$.
Let $Q_{\hbox{\rm dest}}$ be the subset of ${\Q}[t]_{\dim X-1, +}$ of
polynomials arising in that way. Pick a polynomial $\sigma_\infty$ for
which the conclusion of Proposition~\ref{bound30} holds.
%%%
\begin{Lem} The set $Q_{\hbox{\rm dest}}\cap \{\ \sigma\in
{\Q}[t]_{\dim X-1, +}\,|\, \sigma\le\sigma_{\infty}\,\}$
is finite.
\end{Lem}
%%%
\it Proof\rm.
The assumption $\sigma_\K\le\sigma_{\infty}$
provides a lower bound for $\mu(\K)$.
Since the possible coherent sheaves $\E$ vary in a bounded
family, by Proposition~\ref{bound10}, i.e., they are all quotients
of the sheaf $\O_X(-n)^{\oplus v}$, for some large $n$ and $v$, the
lemma
follows from a result of Grothendieck's (\cite{HL2}, Lemma~1.7.9).
$\square$\par\bigskip
%%%
Let $\sigma_1^\p<\cdots<\sigma_t^\p$ be the polynomials
in $Q_{\hbox{\rm dest}}$ which are smaller than $\sigma_{\infty}$. This
gives rise to
"open intervals" $I_0:=\{\, \sigma\, |\,\sigma<\sigma_1^\p\,\}$,
$I_i:=\{\, \sigma\, |\, \sigma_i^\p<\sigma<\sigma_{i+1}^\p\,\}$,
$i=1,...,t-1$, and
$I_t:=\{\, \sigma\, |\, \sigma_t^\p<\sigma\,\}$.
%%%
\begin{Lem}
Let $\sigma_1$ and $\sigma_2$ be two
polynomials in
${\Q}[t]_{\dim X-1, +}\setminus Q_{\hbox{\rm dest}}$.
If $\sigma_1$ and $\sigma_2$ lie both in one of the $I_i$,
then ${\cal {\cal F}{\cal H}}_{P/\G/\H}^{\sigma_1-ss}\cong
{\cal {\cal F}{\cal H}}_{P/\G/\H}^{\sigma_2-ss}$.
\end{Lem}
%%%
\it Proof\rm. By the assumption that $\sigma_1$ and $\sigma_2$
do not lie in $Q_{\hbox{\rm dest}}$, we have
${\cal {\cal F}{\cal H}}_{P/\G/\H}^{\sigma_i-ss}={\cal {\cal F}{\cal H}}_{P/\G/\H}^{\sigma_i-s}$,
$i=1,2$. Therefore, it is enough to show that a framed Hitchin pair
of type $(P,\G,\H)$ is $\sigma_1$-stable if and only if it is
$\sigma_2$-stable. First, let $i<t$. Suppose $\sigma_1<\sigma_2$. If
$(\E,\eps,\phi,\psi)$
is $\sigma_1$-stable but not $\sigma_2$-stable, then there must
be a $\phi$-invariant subsheaf $\K\subset\ker\psi$
with $\sigma_1<\sigma_\K<\sigma_2$. The saturation of $\K$ has the
same property. But since there is no such polynomial, by assumption,
$(\E,\eps,\phi,\psi)$ is also $\sigma_2$-stable.
Next, let $(\E,\eps,\phi,\psi)$ be $\sigma_2$-stable.
If it was not $\sigma_1$-stable, there would be a saturated
$\phi$-invariant
subsheaf $\F$ which is not contained in $\ker\psi$ which
$\sigma_1$-destabilizes
$(\E,\eps,\phi,\psi)$.
Define $\sigma_\F\in {\Q}[t]_{\dim X-1, +}$ by the condition
$$
{P_\F\over\rk\F} -{\sigma_\F\over\rk\F}\q =\q
{P_\E\over\rk\E} -{\sigma_\F\over\rk\E}.
$$
Then $\sigma_1<\sigma_\F<\sigma_2$. If we choose $\F$ such that
$\sigma_\F$ becomes maximal, then
$(\E,\eps,\phi,\psi)$ will be properly
$\sigma_\F$-semistable.
Its associated graded object possesses a $\overline{\phi}$-invariant
saturated subsheaf $\K\subset\ker\overline{\psi}$ with
$\sigma_\K=\sigma_\F$. This is again an impossibility.
In the remainning case $i=t$, we may assume that either $\sigma_1$
or $\sigma_2$ agrees with $\sigma_\infty$. In the former case, i.e.,
$\sigma_\infty<\sigma_2$,
the assertion follows from Proposition~\ref{bound30}.
If $\sigma_1<\sigma_\infty$, then the same argumentation as before can
be applied.
$\square$\par\bigskip
%%%
Now, pick for each $i\in\{\,0,...,t\,\}$ a polynomial $\sigma_i\in I_i$.
Observe that every $(\E,\eps,\phi,\psi)$ which is $\sigma_i$-stable
is also $\sigma_{i-1}^\p$- and $\sigma_i^\p$-semistable.
Therefore, we obtain a diagram
$$
\begin{array}{cccccccccl}
    {\cal {\cal F}{\cal H}}^{\sigma_0-s}_{P/\G/\H}& &\hskip 1cm
& & {\cal {\cal F}{\cal H}}^{\sigma_{t}-s}_{P/\G/\H} \\
  \swarrow \qquad \q \searrow &\swarrow&\dotfill & \searrow &
\swarrow \qquad \qquad  \searrow\id\\
{\cal H}^{ss}_{P/\G} \qquad {\cal {\cal F}{\cal H}}^{\sigma_1^\p-ss}_{P/\G/\H} & &&&
{\cal {\cal F}{\cal H}}^{\sigma_{t}^{\p}-ss}_{P/\G/\H} \qquad
{\cal {\cal F}{\cal H}}_{P/\G/\H}^{\sigma_t-s}
.\\
\end{array}
$$
Here, ${\cal H}_{P/\G}^{ss}$ is the moduli space of semistable Hitchin pairs.
%%%
\subsection*{\q\sc The parameter space and the group actions}
%%%
To avoid further technicalities, we will from now on assume that $X$
is a curve.
%%%
\begin{Ass}
\label{Ass10}
For any $n\ge n_2$, and any semistable oriented framed Hitchin pair
$(E,\eps,\delta,\phi,\psi)$ of type $(d,r,\G,\H,\n)$:
\begin{itemize}
\item $\H(n)$ is globally generated.
\item $E(n)$ is globally generated, and the first cohomology group
      of $E(n)$ vanishes.
\item The conclusion of Proposition~\ref{Sect} holds for all
      positive rational numbers $\sigma$, and
      also the analogous assertion for semistable
      Hitchin pairs of type $(d,r,\G)$.
\item Fix a positive rational number $\sigma_{\infty}$, for which the
      conclusion of Proposition~\ref{bound30} holds. Then, 
      $[d+r(n+1-g)]/2\ge \sigma_\infty$.
\end{itemize}
\end{Ass}
%%%
Again, $n_2=0$ is assumed. We also adopt~\ref{Ass2}.
We start as in the construction of the parameter space for
framed Hitchin pairs.
As explained in the Preliminaries, the universal quotient
${\frak q}_{\frak P}\colon V\otimes\O_{{\frak P}\times X}\lra
{\frak E}_{\frak P}$ defines a
morphism $d({\frak E}_{\frak P})\colon
{\frak P}\lra\Pic X$. By the universal property
of the Picard scheme, 
${\cal D}_{\frak P}:=\pi_{{\frak P}*}\bigl(
\underline{\Hom}(\det({\frak E}_{\frak P}),
\n[{\frak E}_{\frak P}])\bigr)$ is
invertible, and for any point ${\bf p}\in {\frak P}$,
$$
{\cal D}_{\frak P}\langle{\bf p}\rangle\q\cong\q\Hom\bigl(
\det({\frak E}_{{\frak P}|\{{\bf  p}\}\times X}),\n[
{\frak E}_{{\frak P}|\{{\bf  p}\}\times X}]\bigr).
$$
Set $\widehat{\frak T}:=
\P\bigl({\cal D}_{\frak P}^\vee\oplus
S^r\Hom(V,H^0(\H))^\vee\otimes\O_{\frak P}\bigr)$.
Then, we can construct our parameter space as a
closed subscheme ${\frak T}$ of
$\widehat{\frak T}$. Note that,
outside the closed subscheme $\P({\cal D}_{\frak P}^\vee)$
of  $\widehat{\frak T}$, we can
extract the $r$-th root of  the tautological line bundle.
From this, it is clear that there is
a universal family $({\frak E}_{\frak T},\eps_{\frak T},
\delta_{\frak T},\phi_{\frak T},
\widehat{\psi}_{\frak T},{\frak M}_{\frak T},{\frak N}_{\frak T})$
on ${\frak T}\times X$ which has the local universal property.
${\frak T}_0$ is the open subscheme mapping to ${\frak Q}_0$, and
${\frak T}^{(s)s}_0$ is the open subscheme parametrizing
the (semi)stable oriented framed Hitchin pairs.
%%%
\begin{Rem}
Since $X$ curve, the quasi-projective scheme ${\frak Q}_0$ is smooth.
Therefore, the restriction of the universal quotient to ${\frak Q}_0
\times X$ is locally free of rank $r$.
\end{Rem}
%%% 
There is a natural right action by $\SL(V)$ on
${\frak T}$, and the universal family on ${\frak T}\times X$
comes again with an $\SL(V)$-linearization.
Remark~\ref{Aut}~i) and~\ref{Obs200}~iii) show that all the stabilizers
of points in ${\frak T}_0$ which are represented by
stable oriented framed Hitchin pairs are indeed finite.
This is because there are only finitely multiples of
${\id}_V$ in $\SL(V)$.
We must construct the good (geometric) quotient
${\frak T}_0^{(s)s}\catqot\SL(V)$. 
%%%
\subsubsection*{\q The Gieseker map}
%%%
We let ${\cal J}\subset\Pic X$ be 
the Jacobian of degree $d$ line bundles on $X$ and 
$\n_{\cal J}$ be the restriction of $\n$ to ${\cal J}\times X$.
If $d>2g-2$, then ${\frak A}_{\cal J}:=
\pi_{{\cal J}*}\n_{\cal J}$ and ${\frak A}^\p_{\cal J}:=
\pi_{{\cal J}*}(\n_{\cal J}\otimes \pi_X^*\O_X(m))^{\oplus u}$ 
are locally free.
From the universal family, we get, on ${\frak T}_0\times X$,
homomorphisms
%%%
\begin{eqnarray*}
{\frak E}_{{\frak T}_0}\otimes {\frak E}_{{\frak T}_0}^\vee
&\lra &\pi_{{\frak T}_0}^*{\frak N}_{{\frak T}_0}\otimes 
\pi_X^*\O_X(m)^{\oplus u}\\
\O_{{\frak T}_0\times X} &\lra & \pi_{{\frak T}_0}^*{\frak N}_{{\frak T}_0}.
\end{eqnarray*}
%%%
Observe ${\frak E}_{{\frak T}_0}^\vee\cong 
\bigwedge^{r-1}{\frak E}_{{\frak T}_0}\otimes \det
({\frak E}_{{\frak T}_0})^\vee$, because ${\frak E}_{{\frak T}_0}$
is locally free.
Using the surjection $V\otimes \O_{{\frak T}_0\times X}\lra 
{\frak E}_{{\frak T}_0}$, we obtain
%%%
\begin{eqnarray*}
V\otimes \bigwedge^{r-1} V\otimes\O_{{\frak T}_0\times X}
&\lra &\det({\frak E}_{{\frak T}_0})\otimes  
\pi_{{\frak T}_0}^*{\frak N}_{{\frak T}_0}\otimes \pi_X^*\O_X(m)^{\oplus u}\\
\bigwedge^rV\otimes \pi_X^*\O_X(m)^{\oplus u}
&\lra &
\det({\frak E}_{{\frak T}_0})\otimes
\pi_{{\frak T}_0}^*{\frak N}_{{\frak T}_0}\otimes \pi_X^*\O_X(m)^{\oplus u}.
\end{eqnarray*}
Now, project all this to ${\frak T}_0$, so that you get
%%%
\begin{eqnarray*}
V\otimes \bigwedge^{r-1} V\otimes\O_{{\frak T}_0}
&\lra&
{\frak A}_{{\frak T}_0}^\p\otimes {\frak N}_{{\frak T}_0}
\otimes \L_{{\frak T}_0}\\
\bigwedge^r V \otimes M\otimes\O_{{\frak T}_0} 
&\lra & {\frak A}_{{\frak T}_0}^\p\otimes{\frak N}_{{\frak T}_0}
\otimes \L_{{\frak T}_0}.
\end{eqnarray*}
%%%
Here, ${\frak A}^\p_{{\frak T}_0}$ is the pullback of ${\frak A}^\p_{\cal J}$
under the map $d({\frak E}_{{\frak T}_0})$, and
$\L_{{\frak T}_0}$ is some linearized line bundle.
These data define an $\SL(V)$-equivariant morphism
$$
{\bf t}_1\colon {\frak T}_0\lra {\bf P}_1:=
\P\bigl(\Hom(\bigwedge^rV\otimes M\otimes\O_{\cal J}\oplus V
\otimes \bigwedge^{r-1}V\otimes 
\O_{\cal J},{\frak A}_{\cal J}^\p)^\vee\bigr)
$$
which factorizes over an injective morphism ${\frak P}_0\lra {\bf P}_1$,
and ${\bf t}_1^*\O_{{\bf P}_1}(1)={\frak N}_{{\frak T}_0}\otimes
\L_{{\frak T}_0}$. Next, we have a look at the data defined by the
orientation and the framing, i.e., at
$$
\begin{array}{rcccl}
\bigwedge^rV\otimes\O_{{\frak T}_0} &\lra &
\det({\frak E}_{{\frak T}_0}) &\lra & \n[{\frak E}_{{\frak T}_0}]\otimes
\pi_{{\frak T}_0}^*{\frak M}_{{\frak T}_0}\phantom{.}\\
S^r V\otimes \O_{{\frak T}_0} &\lra &
S^r {\frak E}_{{\frak T}_0} &\lra & \pi_X^* S^r \H \otimes
\pi_{{\frak T}_0}^*{\frak M}_{{\frak T}_0}.\\
\end{array}
$$
%%%
Projecting these to ${{\frak T}_0}$, provides us
with
%%%
\begin{eqnarray*}
\bigwedge^r V\otimes \O_{{\frak T}_0} &\lra & {\frak A}_{{\frak T}_0}
\otimes {\frak M}_{{\frak T}_0}\\
S^r\Hom(V, H^0(\H))\otimes \O_{{\frak T}_0} &\lra & {\frak M}_{{\frak T}_0},
\end{eqnarray*}
%%%
and, thus, with a morphism
%%%
$$
{\bf t}_2\colon {\frak T}_0 \lra {\bf P}_2:=
\P\bigl(\Hom(\bigwedge^r V\otimes \O_{\cal J}, {\frak A}_{\cal J})^\vee
\oplus S^r\Hom(V, H^0(\H))^\vee\otimes \O_{\cal J}\bigr),
$$
such that ${\bf t}_2^*\O_{{\bf P}_2}(1)={\frak M}_{{\frak T}_0}$.
The resulting $\SL(V)$-equivariant and injective homomorphism
$$
{\bf t}\colon {\frak T}_0
\lra
{\bf P}_1\times_{\cal J}{\bf P_2}\lra {\Bbb T}:={\bf P}_1\times {\bf P}_2
$$
is our Gieseker map.
We linearize the $\SL(V)$-action on ${\Bbb T}$ in a very ample line 
bundle of the form $\O_{\Bbb T}(1,1)\otimes$(pullback
of a very ample line bundle on ${\cal J}$). 
%%%
\subsubsection*{\q The semistable points in ${\Bbb T}$}
%%%
The key step to the construction of the moduli spaces is
%%%
\begin{Thm}
\label{Stab100}
Let $t=([q\colon V\otimes \O_X\lra E],[\eps,\phi],[\delta,\psi])$
be a point in ${\frak T}_0$. Then the associated point 
${\bf t}(r)\in {\Bbb T}$ is (semi/poly)stable w.r.t.\
the given linearization if and only
$(E,\eps,\delta,\phi,\psi)$ is a (semi/poly)stable oriented framed
Hitchin pair of type $(d,r,\G,\H,\n)$.
\end{Thm}
%%%
First, observe that the action on ${\cal J}$ is trivial,
so that ${\bf t}(t)$ will be (semi)stable
if and only it is (semi)stable 
in ${\Bbb T}_t:={\bf P}_{1,t}\times{\bf P}_{2,t}$
w.r.t.\ the linearization in $\O(1,1)$
where
\begin{eqnarray*}
{\bf P}_{1,t} &:=&
\P\bigl(\Hom(\bigwedge^rV\otimes M\oplus V\otimes \bigwedge^{r-1}V,
H^0(\det(E)(m))^{\oplus u})^\vee\bigr);
\\
{\bf P}_{2,t} &:=&
\P\bigl(\Hom(\bigwedge^rV,H^0(\det(E)))^\vee \oplus
S^r\Hom(V,H^0(\H))^\vee\bigr).
\end{eqnarray*}
%%%
Next, we introduce on ${\bf P}_{2,t}$ the $\C^*$-action
which multiplies the second component by $z$.
Then, there is a family of linearizations of the
$\C^*$-action in $\O_{{\bf P}_{2,t}}$ parametrized by
natural numbers $e,k$ with $0\le e\le k$ \cite{OST}.
We look at the quotients of ${\bf P}_{1,t}\times {\bf P}_{2,t}$
by these linearized $\C^*$-actions.
If $e=0$, then the quotient is
${\bf T}^0_t:=
{\bf P}_{1,t}\times \P\bigl(\Hom(\bigwedge^rV,H^0(\det(E)))^\vee\bigr)$
with induced polarization $\O(1,1)$.
If $e=k$, then the quotient is 
$$
{\bf T}^\infty_t:=
{\bf P}_{1,t}\times \P\bigl(S^r\Hom(V,H^0(\H))^\vee\bigr)
$$
with induced polarization $\O(1,1)$.
In the other cases, the quotient is
$$
{\bf T}_t:={\bf P}_{1,t}\times 
\P\bigl(\Hom(\bigwedge^rV,H^0(\det(E)))^\vee\bigr)
\times \P\bigl(S^r\Hom(V,H^0(\H))^\vee\bigr)
$$
with induced polarization
$
L_k^e:=\O(k,k-e,e).
$
Define 
$$
\sigma_k^e:= {p\over 2}\cdot {e\over k}. 
$$
Note that, by~\ref{Ass10}, for a given (positive) $\sigma<\sigma_\infty$, 
we can find $0<e<k$ satisfying $\sigma=(p/2)(e/k)$.
By the Preliminaries and Remark~\ref{Obs200}~iii), 
Theorem~\ref{Stab100} now reduces to the following
%%%
\begin{Thm}
\label{Stab200}
{\rm i)} The associated point in
${\bf T}^0_{t}$
is (semi/poly)stable if and only
$(E,\eps,\phi)$ is a (semi/poly)stable Hitchin pair of
type $(d,r,\G)$.
\par
{\rm ii)} The associated by point in
$
{\bf T}^{\infty}_{t}$ 
is (semi/poly)stable
if and only if there are no $\phi$-invariant subbundles $F$
of $E$ which are contained in the kernel of $\psi$.
\par
{\rm iii)}
The associated point in
$
{\bf T}_{t}
$
is (semi/poly)stable w.r.t.\ the linearization $L_k^e$
if and only if $(E,\eps,\phi,\psi)$ is a $\sigma_k^e$-(semi/poly)stable
framed Hitchin pair of type $(d,r,\G,\H)$.
\end{Thm}
%%%
\subsubsection*{\q The $\P\bigl(S^r\Hom(V,H^0(\H))^\vee\bigr)$-component}
%%%
By definition, the image of the point ${\bf t}(t)$ in that
space lies in the image of ${\bf R}$ under the $r$-th
Veronese map.
Therefore, the weights are those in ${\bf R}$ multplied by $r$.
%%%
\subsubsection*
{The $\P\bigl(\Hom(\bigwedge^rV,H^0(\det(E))^\vee\bigr)$-component}
%%%
Let $v_1,...,v_p$ be a basis for $V$.
Set $E_i:=q(\langle\,v_1,...,v_i\,\rangle\otimes\O_X)$.
For a one parameter subgroup $\la$ of $\SL(V)$ which is given
w.r.t.\ that basis by weights $\gamma_1\le\cdots\le\gamma_p$,
one computes
$
\mu([q],\la)= -\sum_{i=1}^{p}(\rk E_i-\rk E_{i-1})\gamma_i,
$
in particular, $\mu([q],\la^{(i)})=p\rk E_i-ir$.
%%%
\subsubsection*{\q The ${\bf P}_{1,t}$-component}
%%%
Fix basess $v_1,...,v_p$ of $V$ and  
$m_1,...,\allowbreak m_\mu$ of $H^0(\O_X(m))$,
let $m_1^1,...,\allowbreak m_\mu^1,...,m_1^u,...,m_\mu^u$
be the resulting basis for $M$,
and $\ev\colon \bigwedge^r V\otimes M\lra H^0(\det(E)(m))$ is the natural
map.
Let's look at some special elements in the space $\Hom(\bigwedge^rV\otimes M\oplus
V\otimes\bigwedge^{r-1}V, H^0(\det(E)(m)^{\oplus u}))$.
For each ordered set $I$ of $r$ elements 
in $\{\,1,...,p\,\}$, each $k\in \{\,1,...,\mu\,\}$,
\and $l\in\{\, 1,...,u\,\}$, we define $S_{I,k,l}$
as the element which maps
$(v_{\iota_1}\wedge\cdots\wedge v_{\iota_r})\otimes m^l_k$
to $\ev((v_{\iota_1}\wedge\cdots\wedge v_{\iota_r})\otimes m^l_k)$
and is zero on all other basis elements of $\bigwedge^rV\otimes M$
and also zero on $V\otimes \bigwedge^{r-1}V$.
This element is an eigenvector for the action of the maximal torus
defined by $v_1,...,v_p$.
Indeed, if $\la$ is given by weights $\gamma_1,...,\gamma_p$,
then it acts on $S_{I,k,l}$ with weight
$\gamma_{\iota_1}+\cdots+\gamma_{\iota_r}$.
In the same way, for $I,k,l$ as before and $i,j\in\{\,1,...,r\,\}$,
we define $\Theta_{I,k,l}^{ij}$ as the element which
maps $v_{\iota_i}\otimes (v_{\iota_1}\wedge\cdots\wedge\hat{v}_{\iota_j}
\wedge\cdots\wedge v_{\iota_r})\otimes m_k^l$
to 
$\ev((v_{\iota_1}\wedge\cdots\wedge
v_{\iota_r})\otimes m_k^l)$ and is zero on $\bigwedge^r V\otimes M$ and
all other basis vectors of $V\otimes \bigwedge^{r-1}V$.
These are also eigenvectors, and $\la$ as above acts with weight
$\gamma_{\iota_1}+\cdots+\gamma_{\iota_r}-
\gamma_{\iota_j}+{\gamma}_{\iota_i}$.
By definition, for any $t\in{\frak T}_0$, the component ${\bf t}_1(t)$
lies in the linear subspace of ${\bf P}_{1,t}$ which is spanned
by the $S_{I,k,l}$ and $\Theta^{ij}_{I,k,l}$, and thus the computation
of weights is analogous to that in the first part of this paper.
\par
After these preparations, it is clear that i) and iii) in 
Theorem~\ref{Stab200} can be 
proved in exactly the same way as Theorem~\ref{SemStab} in the first
part of this paper.
In order to see also ii), we first observe that computations in
the space ${\bf T}^\infty_t$
give that the point ${\bf t}_1(t)$ is (semi/poly)stable 
if and only if $(E,\eps,\phi,\psi)$ is $\sigma^*$-(semi/poly)stable
where $\sigma^*=p/2$.
By Assumption~\ref{Ass10}, 
$\sigma^*\ge\sigma_{\infty}$, so that we can conclude
by Proposition~\ref{bound30}. 
Finally, a standard argument shows
%%%
\begin{Prop}
The map ${\bf t}_{{\frak T}_0^{ss}}\colon {\frak T}_0^{ss}
\lra {\Bbb T}^{ss}$ is a finite morphism.
\end{Prop}
%%%
\subsubsection*{\q The outcome}
%%%
The summary of the results of the previous paragraphs is given by
%%%
\begin{Thm}
The good quotient ${\frak T}_0^{ss}\catqot\SL(V)$ exists.
It is a
projective scheme, and the open subscheme
${\frak T}_0^s\catqot\SL(V)$ is a geometric quotient.
\end{Thm}
%%%
\subsection*{\q\sc The moduli spaces}
%%%
Let $\mathop{\rm OFH}^{(s)s}_{d/r/\G/\H/\n}$ be the functor which
assigns
to each noetherian scheme $S$ the set of equivalence classes of families
of (semi)stable oriented framed Hitchin pairs of type $(d,r,\G,\H,\n)$
which are parametrized by $S$, and define the
closed
subfunctor 
$\mathop{\rm OFH}^{(s)s}_{d/r/\G/\H/\n/\hbox{\rm symm}}$ 
of (semi)stable symmetric oriented framed Hitchin pairs.
Finally, we set
${\cal O{\cal F}{\cal H}}^{(s)s}_{d/r/\G/\H/\n}:={\frak T}_0^{(s)s}\catqot\SL(V)$.
%%%
\begin{Thm}
{\rm i)} There is a natural transformation $\vartheta^{(s)s}$
of $\mathop{\rm OFH}^{(s)s}_{d/r/\G/\H/\n}$ into the functor
of points of ${\cal O{\cal F}{\cal H}}^{(s)s}_{d/r/\G/\H/\n}$ which is minimal in the
usual sense (see {\rm i)} of Thm.~\ref{Main}), so that
${\cal O{\cal F}{\cal H}}^{s}_{d/r/\G/\H/\n}$ is a coarse moduli scheme for
stable oriented framed Hitchin pairs of type $(d,r,\G,\H,\n)$. The map
$\vartheta^{ss}(\C)$ induces a bijection between the set of
$S$-equivalence classes of semistable
oriented framed Hitchin pairs of type $(d,r,\G,\H,\n)$
and the set of closed points of ${\cal O{\cal F}{\cal H}}^{ss}_{d/r/\G/\H/\n}$.
There is also a proper {\rm generalized Hitchin map}
$$
\widehat{\xi}\colon {\cal O{\cal F}{\cal H}}^{ss}_{d/r/\G/\H/\n}
\lra \widehat{\Xi}_\G.
$$
\par
{\rm ii)}
There is a closed subscheme ${\cal O{\cal F}{\cal H}}^{(s)s}_{d/r/\G/\H/\n/\hbox{\rm
symm}}$ of ${\cal O{\cal F}{\cal H}}^{(s)s}_{d/r/\G/\H/\n}$ such that the analogues
to {\rm i)} w.r.t.\ the functor
$\mathop{\rm OFH}^{(s)s}_{d/r/\G/\H/\n/\hbox{\rm symm}}$ hold true.
Furthermore, there is a {\rm Hitchin map}
$$
\widehat{\chi}\colon
{\cal O{\cal F}{\cal H}}^{ss}_{d/r/\G/\H/\n/\hbox{\rm
symm}}\lra \widehat{\P}_\G,
$$
mapping a closed point of the moduli space to the characteristic
polynomial of a representing oriented framed Hitchin pair.
The Hitchin map clearly is proper.
\end{Thm}
%%%
\begin{Ex}
We return to the setting of Example~\ref{Ex1}.
Let ${\cal O{\cal F}{\cal H}}$ be the master space. This time,
we have to determine the
quotient $\P(M_1^\vee\oplus M_2^\vee\oplus \C)\catqot (\C^*\times\SL(V))$.
Denote the coordinates by
$(l_0,l_{1,1},l_{1,2},l_{2,1},l_{2,2},\allowbreak s_1,s_2,s_3)$.
The $\SL(V)$-nullforms are cut out by the equations
$H_i=0$, $i=0,...,4$, where $H_0$, $H_1$, $H_2$, and $H_3$
are as before, and $H_4:=s_3$.
Set $g_0:= H_0^2$, $g_1:= H_1^2$, and $g_2:=H_2$.
It follows easily that the master space ${\cal O{\cal F}{\cal H}}$ is
isomorphic to $\P_2\times\P_1$ with coordinates
$([g_0:g_1:g_2], [H_3:H_4])$.
The Hitchin space ${\cal H}$ is the $\C^*$-quotient of the open subset
$H_4\neq 0$, and the space ${\cal {\cal F}{\cal H}}$ is the
$\C^*$-quotient of the open subset
$H_3\neq 0$.
\end{Ex}
%%%
\subsection*{\q\sc The $(\C^*\times\C^*)$-action on the moduli space}
%%%
On 
%our master space 
${\cal O{\cal F}{\cal H}}:={\cal O{\cal F}{\cal H}}^{ss}_{d/r/L/\H/{\cal N}}$,
there are two $\C^*$-actions which commute with each other:
First, there is the $\C^*$-action which comes from multiplying
the twisted endomorphism by a scalar factor. Second, we can multiply
the framing by a non-zero complex number, this yields the
second $\C^*$-action.
As explained in \cite{OST}, it is important to study the fixed
point sets of those $\C^*$-actions, the so-called \it
varieties of reductions\rm .
The fixed point set of the first $\C^*$-action contains two obvious
components.
The first one is ${\cal M}^{ss}_{d/r/\H/{\cal N}}$, the master
space of semistable oriented framed bundles as constructed in~\cite{OST},
corresponding to the points with $\phi=0$.
The second one is
${\cal O{\cal F}{\cal H}}_{\infty}:={\cal O{\cal F}{\cal H}}_{\neq 0}\catqot\C^*$, where
${\cal O{\cal F}{\cal H}}_{\neq 0}$ is the open subset where $\eps\neq 0$.
\par
The fixed point set of the second $\C^*$-action looks as follows:
First, there is ${\cal H}^{ss}_{d/r/\G}=\{\, \psi=0\,\}$,
the moduli space of semistable Hitchin pairs
of type $(d,r,\G)$. Second, there is ${\cal {\cal F}{\cal H}}_{d/r/L/\H}^{\sigma_\infty-ss}$
embedded as the part $\{\,\delta=0\,\}$.
Third, there is the set of the stable points of the form
$((E^\p,\eps,\phi^\p,0)\oplus
(E^{\p\p},\eps,\phi^{\p\p},\psi),\delta)$ which has $t$
components.
From the GIT-process, these $\C^*$-actions come with 
natural linearizations in an ample line bundle on ${\cal O{\cal F}{\cal H}}$, 
let $l$ be the one
of the second $\C^*$-action, then,  as in \cite{OST}, 
one can now conclude
%%%
\begin{Thm} 
For $k>0$ and $e\in\Z$,
let $l_k^e$ be the modification of the linearization $l$ as described
in Part~I of \cite{OST}. Then the GIT-quotients
${\cal O{\cal F}{\cal H}}\catqot_{l_k^e}\C^*$ run through the
moduli spaces ${\H}^{ss}_{d/r/\G}$, and ${\cal {\cal F}{\cal H}}_{d/r/\G/\H}^{\sigma-ss}$,
$\sigma\in{\Q}_{>0}$. In particular, the chain of flips described
in a previous section is a chain of $\C^*$-flips.
\end{Thm}
%%%
\subsubsection*{\q Acknowledgements}
%%%
Starting point of this paper was the suggestion of Professor Okonek to
study the objects of Stupariu's thesis from the algebraic viewpoint.
The final form of the paper was worked out during my stays in Barcelona
and Ramat Gan.
The author acknowledges support by AGE --- Algebraic Geometry in
Europe,
contract No.\ ER-BCHRXCT 940557 (BBW 93.0187),
by SNF, Nr.\ 2000 -- 045209.95/1, by a grant of the
Generalitat de Catalunya (\#1996SGR00060), 
and by a grant of the Emmy Noether Institute.
%%%

\it
\par
\bigskip
\begin{flushright}
Bar-Ilan University\par
Department of Mathematics\\ and Computer Science\par
Ramat-Gan, 52900\par
Israel\par
\medskip
and\par
\medskip
Universit\"at GH Essen\par
FB 6 Mathematik und Informatik\par
D-45117 Essen\par
Deutschland\par
\smallskip
\it E-mail address: \tt schmitt1@@cs.biu.ac.il
\end{flushright}
\end{document}